\RequirePackage{fix-cm}
\documentclass[smallextended,envcountsect]{svjour3}
\smartqed

\usepackage{graphicx}
\usepackage{amsmath,amssymb}
\usepackage{array}
\usepackage{enumitem}
\usepackage{xcolor}
\usepackage{bm}
\usepackage{multirow}
\usepackage{mathdots}

\providecommand{\ie}		{\emph{i.e\@.}\xspace}
\providecommand{\eg}		{\emph{e.g\@.}\xspace}

\providecommand{\myurl}[1][]	{\texttt{web.eecs.umich.edu/$\sim$fessler#1}\xspace}

\providecommand{\onweb}[1]	{Available from \myurl.}

\long\def\comment#1{}

\providecommand{\bcent}		{\begin{center}}
\providecommand{\ecent}		{\end{center}}
\providecommand{\benum}		{\begin{enumerate}}
\providecommand{\eenum}		{\end{enumerate}}
\providecommand{\bitem}		{\begin{itemize}}
\providecommand{\eitem}		{\end{itemize}}
\providecommand{\bvers}		{\begin{verse}}
\providecommand{\evers}		{\end{verse}}
\providecommand{\btab}		{\begin{tabbing}}	
\providecommand{\etab}		{\end{tabbing}}

\newcounter{blist}
\providecommand{\blistmark}	{\makebox[0pt]{$\bullet$}}
\providecommand{\blistitemsep}	{0pt}
\providecommand{\blist}[1][]	{%
\begin{list}{\blistmark}{%
\usecounter{blist}%
\setlength{\itemsep}{\blistitemsep}%
\setlength{\parsep}{0pt}%
\setlength{\parskip}{0pt}%
\setlength{\partopsep}{0pt}%
\setlength{\topsep}{0pt}%
\setlength{\leftmargin}{1.2em}%
\setlength{\labelsep}{0.5\leftmargin}
\setlength{\labelwidth}{0em}%
#1}
}
\providecommand{\elist}		{\end{list}}

\providecommand{\blistitemsep}	{0pt}
\providecommand{\bjfenum}[1][]	{%
\begin{list}{\bcolor{\arabic{blist}.} }{%
\usecounter{blist}%
\setlength{\itemsep}{\blistitemsep}%
\setlength{\parsep}{0pt}%
\setlength{\parskip}{0pt}%
\setlength{\partopsep}{0pt}%
\setlength{\topsep}{0pt}%
\setlength{\leftmargin}{0.0em}%
\setlength{\labelsep}{1.0\leftmargin}
\setlength{\labelwidth}{0pt}%
#1}
}

\newcounter{blistAlph}
\providecommand{\blistAlph}[1][]
{\begin{list}{\makebox[0pt][l]{\Alph{blistAlph}.}}{%
\usecounter{blistAlph}%
\setlength{\itemsep}{0pt}\setlength{\parsep}{0pt}%
\setlength{\parskip}{0pt}\setlength{\partopsep}{0pt}%
\setlength{\topsep}{0pt}%
\setlength{\leftmargin}{1.2em}%
\setlength{\labelsep}{1.0\leftmargin}
\setlength{\labelwidth}{0.0\leftmargin}#1}%
}

\newcounter{blistRoman}
\providecommand{\blistRoman}[1][]
{\begin{list}{\Roman{blistRoman}.}{%
\usecounter{blistRoman}%
\setlength{\itemsep}{0.5em}\setlength{\parsep}{0pt}%
\setlength{\parskip}{0pt}\setlength{\partopsep}{0pt}%
\setlength{\topsep}{0pt}%
\setlength{\leftmargin}{4em}%
\setlength{\labelsep}{0.4\leftmargin}
\setlength{\labelwidth}{0.6\leftmargin}#1}%
}

\usepackage{bbm} 
\providecommand{\jfbbm}[1]	{\xmath{\mathbbm{#1}}} 
\providecommand{\qed}[1][0pt]	{\hfill\raisebox{#1}{\inmath{\Box}}} 

\providecommand{\reals}		{\jfbbm{R}}

\providecommand{\inprod}[2]	{\xmath{\mathop{\langle #1,\, #2 \rangle}\nolimits}}

\let\equivsave\equiv
\def\equiv{\xmath{\equivsave}}

\providecommand{\ba}[1]		{\left[ \begin{array}{#1}}
\providecommand{\ea}		{\end{array} \right]}
\providecommand{\be}		{\begin{equation}}
\providecommand{\ee}[1]		{\label{#1}\end{equation}}
\providecommand{\bea}		{\begin{eqnarray}}
\providecommand{\eea}[1]	{\label{#1}\end{eqnarray}}
\providecommand{\beas}		{\begin{eqnarray*}}
\providecommand{\eeas}		{\end{eqnarray*}}
\providecommand{\beals}[1][1]	{\begin{alignat*}{#1}}	
\providecommand{\eeals}		{\end{alignat*}}

\providecommand{\berr}[2]{
\bgroup
\renewcommand{\theequation}{#1}
\be
#2
\ee{e,#1}
\egroup
\ignorespaces
}

\providecommand{\bearr}[2]{
\bgroup
\renewcommand{\theequation}{#1}
\bea
#2
\eea{e,#1}
\egroup
\ignorespaces
}

\providecommand{\inmath}	{\ensuremath}
\providecommand{\xmath}[1]	{\inmath{#1}\xspace}
\providecommand{\bmath}[1]	{\xmath{\bm{#1}}}	

\providecommand{\paren}[1]	{\xmath{\left(#1\right)}}

\providecommand{\xfrac}[2]	{\xmath{\frac{#1}{#2}}}

\providecommand{\pfrac}[3][]	{\xmath{\!\paren{#1\frac{#2}{#3}}}}

\providecommand{\sinc}		{\Ofun{\mathrm{sinc}}}

\usepackage{jf-accent}
\usepackage{jf-fun}
\usepackage{jf-names}

\newcommand{\cred} {} 

\newcommand{\Find} {\operatorname{Find}}

\newcommand{\Symm} {\jfbbm{S}}

\newcommand{\sign} {\operatorname{sign}}
\newcommand{\Soft} {\operatorname{S}}
\newcommand{\gra} {\operatorname{gra}}

\newcommand{\cB} {\mathcal{B}}

\newcommand{\cH} {\mathcal{H}}

\newcommand{\cG} {\mathcal{G}}
\newcommand{\cF} {\mathcal{F}(\cH)}
\newcommand{\cFone} {\mathcal{F}(\cH_1)}
\newcommand{\cFtwo} {\mathcal{F}(\cH_2)}
\newcommand{\cM} {\mathcal{M}(\cH)}
\newcommand{\cMs} {\mathcal{M}_{\mu}(\cH)}
\newcommand{\cC} {\mathcal{C}_{\beta}(\cH)}

\newcommand{\Bd} {\mathcal{B}}

\renewcommand{\P} {\operatorname{P}}

\newcommand{\zero} {\bmath{0}}

\newcommand{\A} {\bmath{A}}
\newcommand{\B} {\bmath{B}}
\newcommand{\C} {\bmath{C}}
\newcommand{\D} {\bmath{D}}

\newcommand{\G} {\bmath{G}}
\renewcommand{\H} {\bmath{H}}
\newcommand{\I} {\bmath{I}}
\newcommand{\J} {\bmath{J}}
\newcommand{\K} {\bmath{K}}
\newcommand{\LL} {\bmath{L}}
\newcommand{\Z} {\bmath{Z}}
\renewcommand{\S} {\bmath{S}}
\renewcommand{\P} {\bmath{P}}

\newcommand{\T} {\bmath{T}}
\newcommand{\M} {\bmath{M}}
\newcommand{\Mlam} {\xmath{\M_\lambda}}

\newcommand{\q} {\bmath{q}}

\renewcommand{\aa} {\bmath{a}}
\newcommand{\bb} {\bmath{b}}
\newcommand{\cc} {\bmath{c}}

\renewcommand{\u} {\bmath{u}}
\newcommand{\vv} {\bmath{v}}

\newcommand{\h} {\bmath{h}}
\newcommand{\x} {\bmath{x}}
\newcommand{\y} {\bmath{y}}
\newcommand{\g} {\bmath{g}}

\newcommand{\z} {\bmath{z}}

\renewcommand{\zet} {\bmath{\zeta}}
\newcommand{\xii} {\bmath{\xi}}
\newcommand{\etaa} {\bmath{\eta}}
\newcommand{\nuu} {\bmath{\nu}}

\usepackage{fullpage}

\usepackage{hyperref}
\hypersetup{
        linktoc=page,
        colorlinks=true,
        linkcolor=blue,
        citecolor=red,
        filecolor=magenta,
        urlcolor=cyan
}

\begin{document}


\title{Accelerated proximal point method 
for maximally monotone operators
\thanks{This work was supported in part by the National Research Foundation of Korea (NRF) grant 
funded by the Korea government (MSIT) (No. 2019R1A5A1028324),
and the POSCO Science Fellowship of POSCO TJ Park Foundation.}
}

\author{
Donghwan Kim
}

\institute{
        Department of Mathematical Sciences, KAIST, Republic of Korea \at
        \email{donghwankim@kaist.ac.kr}
}

\date{Date of current version: \today}

\maketitle

\begin{abstract}
This paper proposes an accelerated proximal point method 
for maximally monotone operators.
The proof is computer-assisted via the performance estimation problem approach.
The proximal point method includes various well-known convex optimization methods,
such as the proximal method of multipliers and the alternating direction method of multipliers,
and thus the proposed acceleration has wide applications.
Numerical experiments are presented to demonstrate the accelerating behaviors.
\end{abstract}

\section{Introduction}

A fundamental tool for finding a root of a monotone operator is the proximal point method
\cite{martinet:70:rdv,rockafellar:76:moa}.
The monotone operator theory is particularly of interest,
since it is closely related to convex functions and convex minimization
\cite{bauschke:11:caa,combettes:18:mot,ryu:16:apo}.
For example, the proximal point method is useful 
when solving ill-conditioned problems or dual problems.
In particular, 
the augmented Lagrangian method
(\ie, the method of multipliers)~\cite{hestenes:69:mag,powell:69:amf}
and the alternating direction method of multipliers (ADMM)~\cite{gabay:76:ada,glowinski:75:slp} 
are instances of the proximal point method
applied to dual problems
\cite{eckstein:88:tlm,eckstein:92:otd,rockafellar:76:ala}.

To improve the efficiency of the proximal point method, 
accelerating its worst-case rate has been of interest
both in theory and in applications
(see \eg,~\cite{alvarez:01:aip,attouch:20:coa,attouch:19:coi,corman:14:agp,golshtein:79:mli,guler:92:npp,lin:18:caf}).
In specific,
inspired by Nesterov's fast gradient method~\cite{nesterov:83:amf,nesterov:88:oaa},
G\"{u}ler~\cite{guler:92:npp} accelerated the worst-case rate of the proximal point method
for convex minimization
with respect to the cost function.
This yields the fast rate $O(1/i^2)$ where $i$ denotes the number of iterations,
compared to the $O(1/i)$ rate of the proximal point method.
However, this acceleration has not been theoretically generalized to the monotone inclusion problem,
and only somewhat empirical accelerations, 
\eg, via the relaxation and the inertia
(\ie, an implicit version of the heavy ball method~\cite{polyak:64:smo},
or equivalently, Nesterov's and G\"{u}ler's accelerating technique~\cite{guler:92:npp,nesterov:83:amf,nesterov:88:oaa})
in~\cite{alvarez:01:aip,attouch:20:coa,attouch:19:coi,corman:14:agp,golshtein:79:mli},
have been studied.
Therefore, this paper studies accelerating the worst-case rate of the proximal point method
with respect to the fixed-point residual
for maximally monotone operators.
This provides the fast $O(1/i^2)$ rate, which improves upon the rate $O(1/i)$ of the proximal point method
\cite{brezis:78:pid,gu:20:tsc}.
The proof is computer-assisted via the performance estimation problem (PEP) approach~\cite{drori:14:pof}
and its extensions~\cite{drori:20:efo,drori:16:aov,gu:19:oto,gu:20:tsc,kim:16:ofo,kim:18:ala,kim:18:gto,kim:20:ote,lieder:20:otc,ryu:20:osp,taylor:19:sfo,taylor:17:ewc,taylor:17:ssc}.

Under the additional strong monotonicity condition,
the proximal point method has a linear rate in terms of the fixed-point residual
\cite{rockafellar:76:moa},
while the proposed acceleration is not guaranteed to have such a linear rate.
Therefore, this paper further employs a restarting technique
(\eg, \cite[Section 11.4]{nemirovski:94:emi}\cite[Section 5.1]{nesterov:13:gmf})
under the strong monotonicity condition.
This has a linear rate, 
and is faster than the proximal point method for some practical cases.

The proposed acceleration of the proximal point method has wide applications.
This provides an acceleration
to the proximal method of multipliers
\cite{rockafellar:76:ala},
the Douglas-Rachford splitting method~\cite{douglas:56:otn,lions:79:saf},
and ADMM
\cite{gabay:76:ada,glowinski:75:slp}.
The proposed result also applies
to a preconditioned proximal point method such as
the primal-dual hybrid gradient (PDHG) 
method~\cite{chambolle:11:afo,chambolle:16:ote,esser:10:agf,he:12:cao},
(\ie, a preconditioned ADMM),
yielding an accelerated PDHG method.
This paper then shows that the proposed acceleration
applies to a forward method for cocoercive operators.
Existing works on accelerating the forward method 
can be found, for example, in~\cite{attouch:19:coa,lorenz:15:aif}.

Section~\ref{sec:prob} reviews maximally monotone operators, the proximal point method and its known accelerations.
Section~\ref{sec:pep,mono} studies the PEP with respect to the fixed-point residual for monotone inclusion problems.
Section~\ref{sec:pep,appm} proposes a new accelerated proximal point method using the PEP.
Section~\ref{sec:restart} considers a restarting technique 
to yield a linear rate, 
under the additional strongly monotone assumption.
Section~\ref{sec:appl} applies the proposed acceleration
to well-known instances of the proximal point method,
such as the proximal method of multipliers,
the PDHG method,
the Douglas-Rachford splitting method,
and ADMM.
Section~\ref{sec:appl} also provides numerical experiments.
Section~\ref{sec:disc} presents that the proposed approach
also accelerates the forward method for cocoercive operators,
and Sect.~\ref{sec:conc} concludes.

\section{Problem and method}
\label{sec:prob}

\subsection{Monotone inclusion problem}

Let $\cH$ be a real Hilbert space equipped with inner product $\inprod{\cdot}{\cdot}$,
and associated norm $||\cdot||$. 
A set-valued operator $\M\;:\;\cH\to2^\cH$ is \emph{monotone} if
\begin{align}
\inprod{\x - \y}{\u - \vv} \ge 0 \text{ for all } (\x,\u),(\y,\vv)\in\gra\M
\label{eq:mono}
,\end{align}
where $\gra\M := \{(\x,\u)\in\cH\times\cH\;:\:\u\in\M\x\}$
denotes the graph of $\M$.
A monotone operator $\M$ is \emph{maximally monotone} if
there exists no monotone operator $\A\;:\;\cH\to2^\cH$ such that $\gra\A$ 
properly contains $\gra\M$.
Let $\cM$ be the class of maximally monotone operators on $\cH$.
In addition,
a set-valued operator $\M\;:\;\cH\to2^\cH$ is \emph{$\mu$-strongly monotone}
for $\mu\in\reals_{++}$,
if
\begin{align}
\inprod{\x - \y}{\u - \vv} \ge \mu||\x - \y||^2
	\text{ for all } (\x,\u),(\y,\vv)\in\gra\M
\label{eq:mono,strongly}
.\end{align}
Let $\cMs$ be the class of maximally and $\mu$-strongly monotone operators on $\cH$.
Also, define
$\cB(\cH,\cG) = \{\LL\;:\;\cH\to\cG\;|\;\LL
\text{ is linear and bounded}\}$
for a real Hilbert space $\cG$ equipped with inner product $\inprod{\cdot}{\cdot}$,
and let $\LL^*\in\cB(\cG,\cH)$ be the adjoint of $\LL\in\cB(\cH,\cG)$ that satisfies
$\inprod{\LL\x}{\y} = \inprod{\x}{\LL^*\y}$
for all $\x\in\cH$ and $\y\in\cG$.

This paper considers the monotone inclusion problem:
\begin{align}
\Find \;\; \x\in\cH \quad \text{subject to} \quad \zero \in \M\x 
\label{eq:MI}
,\end{align}
where $\M\in\cM$ (or $\M\in\cMs$).
This includes convex problems and convex-concave problems;
a subdifferential $\partial f$ of a closed proper convex function $f\;:\;\cH\to\reals\cup\{\infty\}$
is maximally monotone~\cite{minty:64:otm}.
Let $\cF$ be the class of closed proper convex functions on $\cH$.

We assume that the optimal set $X_*(\M):= \{\x\in\cH \;:\; \zero\in\M\x\}$ is nonempty.
We also assume that
the distance between an initial point $\x_0$ and some optimal point $\x_*\in X_*(\M)$ 
is bounded as
\begin{align}
||\x_0 - \x_*|| \le R \quad \text{for a constant } R>0
\label{eq:init,dist}
.\end{align}

\subsection{Proximal point method and its worst-case rates}

Proximal point method was first introduced to convex optimization by Martinet~\cite{martinet:70:rdv},
which is based on the proximal mapping by Moreau~\cite{moreau:65:ped}.
The method was later extended to monotone inclusion problem by Rockafellar~\cite{rockafellar:76:moa}.
The proximal point method for maximally monotone operators
includes
the augmented Lagrangian~\cite{hestenes:69:mag,powell:69:amf},
the proximal method of multipliers~\cite{rockafellar:76:ala},
the Douglas-Rachford splitting method~\cite{douglas:56:otn,lions:79:saf},
and the alternating direction method of multipliers (ADMM)
\cite{gabay:76:ada,glowinski:75:slp},
so studying its worst-case convergence behavior and acceleration is important,
which is of main interest in this paper.

The proximal mapping~\cite{moreau:65:ped} (or the resolvent operator) of an operator $\M$ is defined as
\begin{align}
\J_{\M} := (\I + \M)^{-1}
,\end{align}
where $\I\;:\;\cH\to\cH$ is an identity operator, \ie, $\I(\x)=\x$ for all $\x\in\cH$.
The resolvent operator $\J_{\M}$ is single-valued and firmly nonexpansive 
for $\M\in\cM$~\cite{minty:62:mno}.
The proximal point method~\cite{martinet:70:rdv,rockafellar:76:moa} 
generates a sequence $\{\x_i\}$ by iteratively applying the resolvent operator
with a positive real number $\lambda$ as below.

\fbox{
\begin{minipage}[t]{0.85\linewidth}
\vspace{-10pt}
\begin{flalign*}
&\quad \text{\bf Proximal Point Method} & \\
&\qquad \text{Input: } \M\in\cM,\; \x_0\in\cH,\; \lambda\in\reals_{++}. & \\
&\qquad \text{For } i = 0,1,\ldots & \\
&\qquad \qquad \x_{i+1} = \J_{\lambda\M}(\x_i). & 
\end{flalign*}
\end{minipage}
} \vspace{5pt}

In~\cite[Proposition 8]{brezis:78:pid}, 
the worst-case rate of the proximal point method
with respect to the fixed-point residual 
\begin{align}
||\x - \J_{\lambda\M}(\x)||^2
\label{eq:residual}
\end{align}
was found to satisfy
\begin{align}
||\x_i - \x_{i-1}||^2 \le \frac{R^2}{i}
\label{eq:ppm,mono}
\end{align}
for $i\ge1$.
Very recently in~\cite{gu:20:tsc}, this was improved to
\begin{align}
||\x_i - \x_{i-1}||^2 
\le \left(1 - \frac{1}{i}\right)^{i-1}\frac{R^2}{i}
\label{eq:ppm,rate}
,\end{align}
which is exact when $\dim \cH \ge 2$.
Such exact worst-case with $\dim \cH = 2$ given in~\cite{gu:20:tsc} 
will be visited at the end of Sect.~\ref{sec:pep,appm}.
The bound~\eqref{eq:ppm,rate} is 
asymptotically $e$-times lower than~\eqref{eq:ppm,mono},
where $e$ is Euler's number.
When we additionally assume the $\mu$-strong monotonicity,
the proximal point method
has a linear rate
\cite[Example 23.40]{bauschke:11:caa}
\cite{rockafellar:76:moa}
\begin{align}
||\x_{i+1} - \x_i||^2 \le 
\left(\frac{1}{1+\lambda\mu}\right)^2||\x_i - \x_{i-1}||^2
\label{eq:ppm,rate,strongly}
\end{align}
for $i\ge1$,
which is exact considering the case 
$\M\x = \mu\x$ with $\dim \cH = 1$.

For a convex minimization of $f\in\cF$,
\cite[Conjecture 4.2]{taylor:17:ewc} conjectures
that the proximal point method satisfies
\begin{align}
||\x_i - \x_{i-1}||^2 \le \frac{R^2}{i^2}
\label{eq:ppm,rate,convex}
\end{align}
for $i\ge1$,
which is faster than~\eqref{eq:ppm,rate} for maximally monotone operators.
In addition,
the $O(1/i)$ worst-case rate of the proximal point method 
with respect to the cost function
was studied in~\cite[Theorem 2.1]{guler:91:otc},
and this was improved by a constant $2$
in~\cite[Theorem 4.1]{taylor:17:ewc}
\begin{align}
f(\x_i) - f(\x_*) \le \frac{R^2}{4\lambda i} 
\label{eq:ppm,rate,convex,cost}
\end{align}
for $i\ge1$ and some $\x_* \in X_*(\partial f)$ with $||\x_0 - \x_*|| \le R$.

\begin{remark}
\label{remark:precond}
The results for the proximal point method
can be applied to a \emph{preconditioned} proximal point method.
Let $\LL \in \cB(\cH,\cH)$ be invertible.
Then, $\LL^*\M\LL$ is maximally monotone
for $\M\in\cM$
\cite[Proposition 23.25]{bauschke:11:caa},
and the corresponding proximal point method is
\begin{align}
\tilde{\x}_{i+1} = \J_{\lambda\LL^*\M\LL}(\tilde{\x}_i) = (\I + \lambda\LL^*\M\LL)^{-1}\tilde{\x}_i
\label{eq:ppm,composite}
.\end{align}
Introducing $\x_i = \LL\tilde{\x}_i$ and 
$\P = (\LL\LL^*)^{-1}$ 
yields the following equivalent 
\emph{preconditioned} proximal point method 
\begin{align}
\x_{i+1} = (\P + \lambda\M)^{-1}\P\x_i
\label{eq:ppm,precond}
.\end{align}
So, for example,
the inequality~\eqref{eq:ppm,mono} leads to the \emph{preconditioned} fixed-point residual bound for
the \emph{preconditioned} proximal point method
\begin{align}
\inprod{\P(\x_i - \x_{i-1})}{\x_i - \x_{i-1}} \le \frac{R^2}{i}
\end{align}
for $i\ge1$,
and for some $\x_*\in X_*(\M)$ with $\inprod{\P(\x_0 - \x_*)}{\x_0 - \x_*} \le R^2$.
This is particularly useful
when considering the PDHG method
\cite{chambolle:11:afo,chambolle:16:ote,esser:10:agf}\cite[Lemma 2.2]{he:12:cao},
which is 
an instance of 
a \emph{preconditioned} proximal point method.
We will revisit this in Sect.~\ref{sec:acpm}.
\end{remark}

\subsection{Existing accelerations for proximal point method}

This section reviews existing accelerations of proximal point method for convex minimization
with respect to the cost function.
To the best of our knowledge, there is no other type of proximal point methods
that guarantees accelerated worst-case rates.

For convex minimization,
G\"{u}ler~\cite{guler:92:npp}
developed the following two accelerated versions, inspired 
by Nesterov's fast gradient method~\cite{nesterov:83:amf,nesterov:88:oaa}.
The following is the first accelerated version of the proximal point method in~\cite{guler:92:npp}
which is an instance of FISTA~\cite{beck:09:afi}.
The original version in~\cite{guler:92:npp} includes some variation
with an iteration-dependent $\lambda_i$,
rather than a fixed constant $\lambda$
(see also~\cite{attouch:19:fpm} for choosing $\lambda_i$ appropriate for further acceleration).
This paper focuses on a fixed constant $\lambda$,
and we leave its extension to a varying constant $\lambda_i$
as future work.

\fbox{
\begin{minipage}[t]{0.85\linewidth}
\vspace{-10pt}
\begin{flalign*}
&\quad \text{\bf G\"{u}ler's First Accelerated Proximal Point Method in~\cite[Sec. 2]{guler:92:npp}} & \\
&\qquad \text{Input: } f\in\cF,\; \x_0=\y_0\in\cH,\; \lambda\in\reals_{++},\; t_0=1. & \\
&\qquad \text{For } i = 0,1,\ldots & \\
&\qquad \qquad \x_{i+1} = \J_{\lambda\partial f}(\y_i), & \\
&\qquad \qquad t_{i+1} = \frac{1 + \sqrt{1 + 4t_i^2}}{2}, & \\
&\qquad \qquad \y_{i+1} = \x_{i+1}
                + \frac{t_i-1}{t_{i+1}}(\x_{i+1} - \x_i). &
\end{flalign*}
\end{minipage}
} \vspace{5pt}

The sequence generated by the G\"{u}ler's first accelerated proximal point method 
satisfies~\cite[Theorem 2.3]{guler:92:npp}
\cite[Theorem 4.4]{beck:09:afi}
\begin{align}
f(\x_i) - f(\x_*) \le \frac{R^2}{2\lambda t_{i-1}^2}
	\le \frac{2R^2}{\lambda(i+1)^2}
\label{eq:guler,ppm,rate}
\end{align}
for $i\ge1$
and for some $\x_* \in X_*(\partial f)$ with $||\x_0 - \x_*|| \le R$.
The following is another accelerated proximal point method by G\"{u}ler~\cite{guler:92:npp},
which the formulation is similar to those of the optimized gradient methods
\cite{kim:16:ofo,kim:18:gto,kim:20:ote}. 

\fbox{
\begin{minipage}[t]{0.85\linewidth}
\vspace{-10pt}
\begin{flalign*}
&\quad \text{\bf G\"{u}ler's Second Accelerated Proximal Point Method in~\cite[Appendix]{guler:92:npp}} & \\
&\qquad \text{Input: } f\in\cF,\; \x_0=\y_0\in\cH,\; \lambda\in\reals_{++},\; t_0=1. & \\
&\qquad \text{For } i = 0,1,\ldots & \\
&\qquad \qquad \x_{i+1} = \J_{\lambda\partial f}(\y_i), & \\
&\qquad \qquad t_{i+1} = \frac{1 + \sqrt{1 + 4t_i^2}}{2}, & \\
&\qquad \qquad \y_{i+1} = \x_{i+1}
                + \frac{t_i-1}{t_{i+1}}(\x_{i+1} - \x_i)
                + \frac{t_i}{t_{i+1}}(\x_{i+1} - \y_i). &
\end{flalign*}
\end{minipage}
} \vspace{5pt}

The sequence generated by G\"{u}ler's second accelerated proximal point method 
satisfies~\cite[Theorem 6.1]{guler:92:npp} for $i\ge1$
\begin{align}
f(\x_i) - f(\x_*) \le \frac{R^2}{4\lambda t_{i-1}^2}
	\le \frac{R^2}{\lambda(i+1)^2},
\end{align}
which is twice smaller than~\eqref{eq:guler,ppm,rate}.

\subsection{Main contribution}

To accelerate the worst-case rate of the proximal point method for maximally monotone operators,
the relaxation and the inertia 
(\ie, an implicit version of the heavy ball method~\cite{polyak:64:smo},
or equivalently, Nesterov's and G\"{u}ler's accelerating technique~\cite{guler:92:npp,nesterov:83:amf,nesterov:88:oaa})
have been studied
in~\cite{alvarez:01:aip,attouch:20:coa,attouch:19:coi,corman:14:agp,golshtein:79:mli}.
However, none of them guarantee accelerated rates.
Therefore,
the main contribution of this paper
is to develop a method that has a fast $O(1/i^2)$ rate
with respect to the fixed-point residual,
improving upon the $O(1/i)$ rate of the proximal point method 
in~\eqref{eq:ppm,mono} and~\eqref{eq:ppm,rate}.

This paper considers the following general proximal point method
with step coefficients $\{h_{i+1,k+1}\}_{k=0}^i$
for reusing previous and current updates $\{\x_{k+1} - \y_k\}_{k=0}^i$.
This includes the proximal point method, 
the accelerated methods via the relaxation and the inertia 
\cite{alvarez:01:aip,attouch:20:coa,attouch:19:coi,corman:14:agp,golshtein:79:mli},
and the proposed accelerated method.

\fbox{
\begin{minipage}[t]{0.85\linewidth}
\vspace{-10pt}
\begin{flalign*}
&\quad \text{\bf General Proximal Point Method} & \\
&\qquad \text{Input: } \M\in\cM,\; \y_0\in\cH,\; \lambda\in\reals_{++}. & \\
&\qquad \text{For } i = 0,1,\ldots & \\
&\qquad \qquad \x_{i+1} = \J_{\lambda\M}(\y_i), & \\
&\qquad \qquad \y_{i+1} = \y_i + \sum_{k=0}^i h_{i+1,k+1}(\x_{k+1} - \y_k). & 
\end{flalign*}
\end{minipage}
} \vspace{5pt}

This paper next uses the PEP approach to find the choice of $\{h_{i+1,k+1}\}_{k=0}^i$
that guarantees an accelerated rate.
While the formulation of the general proximal point method 
is inefficient in general,
the proposed accelerated method with the specific choice of $\{h_{i+1,k+1}\}_{k=0}^i$ found by PEP
has an efficient equivalent form.
This form is similar to the other accelerated methods
with the relaxation and/or the inertia
\cite{alvarez:01:aip,attouch:20:coa,attouch:19:coi,corman:14:agp,golshtein:79:mli}.

\section{Performance estimation problem for maximally monotone operators}
\label{sec:pep,mono}

This section uses the performance estimation problem (PEP) approach~\cite{drori:14:pof,taylor:17:ewc,taylor:17:ssc}
to analyze the general proximal point method for maximally monotone operators,
in terms of the fixed-point residual~\eqref{eq:residual}.
This was recently studied in~\cite{gu:20:tsc}
for the proximal point method, providing the exact rate~\eqref{eq:ppm,rate}.
The same authors~\cite{gu:19:oto} also used the PEP to study the exact worst-case rate
for the \emph{ergodic} sequence of the (relaxed) proximal point method for the variational inequalities.
Similarly,~\cite{taylor:17:ewc} used PEP to analyze the worst-case rate 
of the proximal point method for convex minimization
in terms of the fixed-point residual and the cost function,
yielding~\eqref{eq:ppm,rate,convex} and~\eqref{eq:ppm,rate,convex,cost}, respectively.

Building upon~\cite{drori:14:pof,gu:19:oto,gu:20:tsc,taylor:17:ewc,taylor:17:ssc},
the worst-case rate of the general proximal point method after $N$ iterations
for decreasing the fixed-point residual~\eqref{eq:residual} 
under the initial distance condition~\eqref{eq:init,dist} can be computed by
\begin{align}
\max_{\M\in\cM}\max_{\substack{\x_1,\ldots,\x_N\in\cH, \\ 
	\y_0,\ldots,\y_{N-1}\in\cH, \\
	\x_*\in X_*(\M)}}\; & \frac{1}{R^2}||\x_N - \y_{N-1}||^2 
	\label{eq:P} \\
\text{subject to }\; & \x_{i+1} = \J_{\lambda\M}(\y_i), \quad i=0,\ldots,N-1, \nonumber \\
        & \y_{i+1} = \y_i + \sum_{k=0}^i h_{i+1,k+1}(\x_{k+1} - \y_k), \quad i=0,\ldots,N-2, \nonumber \\
        & ||\y_0 - \x_*||^2 \le R^2. \nonumber
\end{align}
This is an infinite-dimensional problem 
due to the constraint $\M\in\cM$,
which is impractical to solve.
PEP in~\cite{drori:14:pof} further introduced a series of steps
that reformulate such impractical problem
into a tractable problem, which we apply to~\eqref{eq:P}
step by step below.

The first step is to reformulate the problem~\eqref{eq:P}
into a finite-dimensional problem.
\cite[Fact 1]{ryu:20:osp} implies that
one can replace $\M\in\cM$ in~\eqref{eq:P}
by a set of inequality constraints~\eqref{eq:mono} for $\M\in\cM$
on the finite number of pairs of points $\{\x_1,\ldots,\x_N,\x_*\}$
without strictly relaxing the problem~\eqref{eq:P}.
In specific, such constraints are
\begin{align}
\inprod{\x_i - \x_j}{\q_i - \q_j} \ge 0,
\end{align}
for all $i,j\in\{1,\ldots,N,*\}$,
with additional variables $\q_i\in\M\x_i$ for $i=1,\ldots,N$
and $\q_* = \zero \in \M\x_*$.
Then the resulting equivalent problem of~\eqref{eq:P} is
\begin{align}
\max_{\substack{\x_1,\ldots,\x_N,\x_*\in\cH, \\ \y_0,\ldots,\y_{N-1}\in\cH, \\
\q_1,\ldots,\q_N\in\cH}}\; & \frac{1}{R^2}||\x_N - \y_{N-1}||^2 \label{eq:Pfin} \\
\text{subject to }\; & \inprod{\x_i - \x_j}{\q_i - \q_j} \ge 0, \quad i<j=1,\ldots,N, \nonumber \\
	& \inprod{\x_i - \x_*}{\q_i} \ge 0, \quad i=1,\ldots,N, \nonumber \\
        & \x_{i+1} = \y_i - \lambda\q_{i+1}, \quad i=0,\ldots,N-1, \nonumber \\
        & \y_{i+1} = \y_i - \lambda\sum_{k=0}^i h_{i+1,k+1}\q_{k+1}, \quad i=0,\ldots,N-2, \nonumber \\
        & ||\y_0 - \x_*||^2 \le R^2. \nonumber
\end{align}
Further removing $\x_i$
and using the change of variables
\begin{align}
\g_i := \frac{\lambda}{R}\q_i, \quad i=1,\ldots,N,
\end{align}
simplify the problem~\eqref{eq:Pfin} as
\begin{align}
\max_{\substack{\y_0,\ldots,\y_{N-1},\x_*\in\cH, \\ \g_1,\ldots,\g_N\in\cH}}\; & ||\g_N||^2 \\
\text{subject to }\; & \frac{1}{R}\inprod{\y_{i-1} - R\g_i - \y_{j-1} + R\g_j}{\g_i - \g_j} \ge 0, \quad i<j=1,\ldots,N, \nonumber \\
	& \frac{1}{R}\inprod{\y_{i-1} - R\g_i - \x_*}{\g_i} \ge 0, \quad i=1,\ldots,N, \nonumber \\
        & \y_{i+1} = \y_i - R\sum_{k=0}^i h_{i+1,k+1}\g_{k+1}, \quad i=0,\ldots,N-2, \nonumber \\
        & ||\y_0 - \x_*||^2 \le R^2. \nonumber
\end{align}

As in~\cite{drori:14:pof,gu:19:oto,gu:20:tsc,taylor:17:ewc,taylor:17:ssc},
we next introduce the Gram matrix
\begin{align}
\Z = \left[\begin{array}{ccccc}
||\g_1||^2 & \inprod{\g_1}{\g_2} & \cdots & \inprod{\g_1}{\g_N} & \frac{1}{R}\inprod{\g_1}{\y_0 - \x_*} \\
\inprod{\g_1}{\g_2} & ||\g_2||^2 & \cdots & \inprod{\g_2}{\g_N} & \frac{1}{R}\inprod{\g_2}{\y_0 - \x_*} \\
\vdots & \vdots& \ddots & \vdots & \vdots \\
\inprod{\g_1}{\g_N} & \cdots & & ||\g_N||^2 & \frac{1}{R}\inprod{\g_N}{\y_0 - \x_*} \\
\frac{1}{R}\inprod{\g_1}{\y_0 - \x_*} & \cdots & & \frac{1}{R}\inprod{\g_N}{\y_0 - \x_*} & \frac{1}{R^2}||\y_0 - \x_*||^2 
\end{array}\right]
\label{eq:Z}
\end{align}
to relax the problem as
\begin{align}
\max_{\Z\in\Symm_+^{N+1}}\; & \tr\{\u_N\u_N^\top\Z\} 
	\label{eq:PZ} \\
\text{subject to }\; & \tr\{\A_{i,j}(\h)\Z\} \le 0, \quad i<j=1,\ldots,N, \nonumber \\
        & \tr\{\B_i(\h)\Z\} \le 0, \quad i=1,\ldots,N, \nonumber \\
        & \tr\{\C\Z\} \le 1, \nonumber
\end{align}
where 
$\{\u_i\}_{i=1}^{N+1}$ is the canonical basis of $\reals^{N+1}$
and
\begin{align*}
\begin{cases}
\A_{i,j}(\h) := (\u_i-\u_j)\odot(\u_i - \u_j)
        - (\u_i - \u_j)\odot\sum_{l=i-1}^{j-2}\sum_{k=0}^lh_{l+1,k+1}\u_{k+1}, 
        & i<j=1,\ldots,N, \\
\B_i(\h) := \u_i\u_i^\top
        - \u_i\odot\u_{N+1}
        + \u_i\odot\sum_{l=0}^{i-2}\sum_{k=0}^lh_{l+1,k+1}\u_{k+1}, 
        & i=1,\ldots,N, \\
\C := \u_{N+1}\u_{N+1}^\top
\end{cases}
\end{align*}
with the outer product operator $\u\odot\vv := \frac{1}{2}(\u\vv^\top + \vv\u^\top)$.
If $\dim\cH \ge N+1$, 
the problems~\eqref{eq:P} and~\eqref{eq:PZ}
are equivalent, based on the following lemma
similar to \cite[Lemma 1]{ryu:20:osp}.
\begin{lemma}
\label{lem:Z}
If $\dim\cH \ge N+1$, then 
\begin{align*}
\Z\in\Symm_+^{N+1} \quad \Leftrightarrow \quad
\exists\; \g_1,\g_2,\ldots,\g_N,\frac{1}{R}(\y_0-\x_*)\in\cH
	\text{ such that } \Z = \text{expression of}~\eqref{eq:Z}
.\end{align*}
\end{lemma}

For simplicity in later analysis,
we discard some constraints as
\begin{align}
\max_{\Z\in\Symm_+^{N+1}}\; & \tr\{\u_N\u_N^\top\Z\} 
	\label{eq:PZrelax} \\
\text{subject to }\; & \tr\{\A_{i-1,i}(\h)\Z\} \le 0, \quad i=2,\ldots,N, \nonumber \\
        & \tr\{\B_N(\h)\Z\} \le 0,  
		\nonumber \\
        & \tr\{\C\Z\} \le 1, \nonumber
\end{align}
which does not affect the result in the paper, \ie,
the optimal values of~\eqref{eq:PZ} and~\eqref{eq:PZrelax}
are found to be numerically equivalent for the method proposed in this paper.
Finally, we construct the associated Lagrangian dual of~\eqref{eq:PZrelax}
\begin{align}
\Bd_{D}(\h) := \qquad 
\min_{a_2,\ldots,a_N,b_N,c\in\reals}\; & c 
	\tag{D} 
	\label{eq:D} \\
\text{subject to }\; & \sum_{i=2}^Na_i\A_{i-1,i}(\h) 
		+ b_N\B_N(\h)
		+ c\C - \u_N\u_N^\top \succeq \zero, \nonumber \\ 
        & a_2,\ldots,a_N,b_N,c \ge 0, \nonumber
\end{align}
where $a_2,\ldots,a_N,b_N,c$ are dual variables 
associated with the constraints of~\eqref{eq:PZrelax}, respectively.
Then, for any given $\h$ for the general proximal point method,
one can compute its (upper bound of) 
worst-case fixed-point residual by numerically solving~\eqref{eq:D} using any SDP solver.
For some choices of $\h$ as for the proximal point method in~\cite{gu:20:tsc}, 
it might be possible to analytically solve~\eqref{eq:D};
\cite{gu:20:tsc} analytically solved~\eqref{eq:D} for the proximal point method
yielding the rate~\eqref{eq:ppm,rate}.
This paper provides another choice of $\h$ that provides an analytical solution to~\eqref{eq:D}
with an accelerated rate.

\section{Accelerating the proximal point method for maximally monotone operators}
\label{sec:pep,appm}

Using the dual problem~\eqref{eq:D},
this section develops an accelerated version of the proximal point method
via PEP:
\begin{align}
\min_{\h} \Bd_{D}(\h) 
\label{eq:HD} \tag{HD}
,\end{align}
which is studied in~\cite{drori:20:efo,drori:14:pof,drori:16:aov,kim:16:ofo,kim:18:ala,kim:18:gto,kim:20:ote}
for certain classes of problems and methods.
The problem is non-convex but convex for the variables $(a_2,\ldots,a_N,b_N,c)$ given $\h$
and for the variables $(c,\h)$ given $(a_2,\ldots,a_N,b_N)$.
Therefore,
we used a variant of alternating minimization 
that alternatively optimizes over $(a_2,\ldots,a_N,b_N,c)$ given $\h$
and over $(c,\h)$ given $(a_2,\ldots,a_N,b_N)$
to find a minimizer
using a SDP solver~\cite{cvxi,gb08}.
Inspired by numerical results,
the following lemma specifies a feasible point of~\eqref{eq:HD} analytically.
We do not have a guarantee that such point is a (unique) minimizer of~\eqref{eq:HD}.

\begin{lemma}
\label{lem:feas}
The following
\begin{align}
h_{i,k} &= \begin{cases}
-\frac{2k}{i(i+1)}, & i=1,\ldots,N-1,\;k=1,\ldots,i-1, \\
\frac{2i}{i+1}, & i=1,\ldots,N-1,\;k=i,
\end{cases}
\label{eq:hik} \\
a_i &= \frac{2(i-1)i}{N^2}, \quad i=2,\ldots,N, \quad 
b_N = \frac{2}{N}, \quad
c = \frac{1}{N^2}
\label{eq:abc}
\end{align}
is a feasible point of~\eqref{eq:D} and~\eqref{eq:HD}.
\begin{proof}
It is obvious that $a_2,\ldots,a_N,b_N,c$ are nonnegative,
so we are only left to show the positive semidefinite condition in~\eqref{eq:D}.
Since
\begingroup
\allowdisplaybreaks
\begin{align*}
&\; \sum_{i=2}^Na_i\A_{i-1,i}(\h) + b_N\B_N(\h) + c\C - \u_N\u_N^\top \\
= &\; \sum_{i=2}^N\frac{2(i-1)i}{N^2}\left[
	(\u_{i-1} - \u_i)\odot(\u_{i-1} - \u_i)
	- (\u_{i-1} - \u_i)\odot\left(\frac{2(i-1)}{i}\u_{i-1}
	- \sum_{k=0}^{i-3}\frac{2(k+1)}{(i-1)i}\u_{k+1}\right)
	\right] \\
	& + \frac{2}{N}\left[\u_N\u_N^\top - \u_N\odot\u_{N+1} 
		+ \u_N\odot\sum_{l=0}^{N-2}\left(\frac{2(l+1)}{l+2}\u_{l+1} 
		- \sum_{k=0}^{l-1}\frac{2(k+1)}{(l+1)(l+2)}\u_{k+1}\right)\right] \\
	& + \frac{1}{N^2}\u_{N+1}\u_{N+1}^\top
	- \u_N\u_N^\top \\
= &\; \sum_{i=2}^{N-1} \left[\frac{2(i-1)i}{N^2} 
	+ \frac{2i(i+1)}{N^2}\paren{1 - \frac{2i}{i+1}}\right]\u_i\u_i^\top 
	+ \left[\frac{2(N-1)N}{N^2} + \frac{2}{N} - 1\right]\u_N\u_N^\top
	+ \frac{1}{N^2}\u_{N+1}\u_{N+1}^\top \\
	& + \sum_{i=2}^{N-1}\left[\frac{2(i-1)i}{N^2}\left(-2 + \frac{2(i-1)}{i}\right) 
		+ \frac{2i(i+1)}{N^2}\frac{2(i-1)}{i(i+1)}\right]\u_{i-1}\odot\u_i \\
	& + \left[\frac{2(N-1)}{N}\paren{-2 + \frac{2(N-1)}{N}}  + \frac{2}{N}\frac{2(N-1)}{N}\right]\u_{N-1}\odot\u_N
		- \frac{2}{N}\u_N\odot\u_{N+1} \\
	& + \sum_{i=3}^{N-1}\sum_{k=0}^{i-3}\left[-\frac{2(i-1)i}{N^2}\frac{2(k+1)}{(i-1)i} + \frac{2i(i+1)}{N^2}\frac{2(k+1)}{i(i+1)}\right]\u_{k+1}\odot\u_i \\
	& + \sum_{k=0}^{N-3}\left[-\frac{2(N-1)N}{N^2}\frac{2(k+1)}{(N-1)N} + \frac{2}{N}\left(\frac{2(k+1)}{k+2}
		- \sum_{l=k+1}^{N-2}\frac{2(k+1)}{(l+1)(l+2)}\right)\right]\u_{k+1}\odot\u_N \\
= &\; \u_N\u_N^\top
	+ \frac{1}{N^2}\u_{N+1}\u_{N+1}^\top 
	- \frac{2}{N}\u_N\odot\u_{N+1} \\
=&\; \left(\u_N - \frac{1}{N}\u_{N+1}\right)\left(\u_N - \frac{1}{N}\u_{N+1}\right)^\top
\succeq \zero,
\end{align*}
\endgroup
the given point is a feasible point of~\eqref{eq:HD}.
\qed
\end{proof}
\end{lemma}

Before providing the worst-case rate of the general proximal point method 
with $\h$ in~\eqref{eq:hik},
we develop its efficient formulation below.
This has a low computational cost per iteration, comparable to that of the proximal point method.
Note that this may not be the only efficient form for $\h$ in~\eqref{eq:hik}.

\fbox{
\begin{minipage}[t]{0.85\linewidth}
\vspace{-10pt}
\begin{flalign*}
&\quad \text{\bf Proposed Accelerated Proximal Point Method for Maximally Monotone Operators} & \\
&\qquad \text{Input: } \M\in\cM,\; \x_0=\y_0=\y_{-1}\in\cH,\; \lambda\in\reals_{++}. & \\
&\qquad \text{For } i = 0,1,\ldots & \\
&\qquad \qquad \x_{i+1} = \J_{\lambda\M}(\y_i), & \\
&\qquad \qquad \y_{i+1} = \x_{i+1} + \frac{i}{i+2}(\x_{i+1} - \x_i) - \frac{i}{i+2}(\x_i - \y_{i-1}). &
\end{flalign*}
\end{minipage}
} \vspace{5pt}

\begin{proposition}
\label{prop:appm}
The sequences $\{\x_i\}$ and $\{\y_i\}$ generated by the general proximal point method
with step coefficients $\{h_{i,k}\}$ in~\eqref{eq:hik}
are identical to the corresponding sequence generated by
the proposed accelerated proximal point method starting from the same initial point.
\begin{proof}
We use induction, and for clarity we use the notation $\x_1',\x_2',\ldots$
and $\y_0',\y_1',\ldots$
for the general proximal point method with~\eqref{eq:hik}.
It is obvious that $\x_0=\y_0'=\y_0$, $\x_1' = \x_1 = \y_1$, and we have
\begin{align*}
\y_1' &= \y_0' + h_{1,1}(\x_1' - \y_0') = \x_1' = \y_1 
.\end{align*}
Similarly, it is obvious that $\x_2' = \x_2$, and we have
\begin{align*}
\y_2' &= \y_1' + \sum_{k=0}^1h_{2,k+1}(\x_{k+1}' - \y_k') 
= \y_1 + \frac{4}{3}(\x_2 - \y_1) - \frac{1}{3}(\x_1 - \y_0) \\
&= \x_2 + \frac{1}{3}(\x_2 - \x_1) - \frac{1}{3}(\x_1 - \y_0) = \y_2.
\end{align*}
It is then also obvious that $\x_3' = \x_3$.
Assuming $\x_l' = \x_l$ for $l=1,\ldots,i+1$ and $\y_l' = \y_l$ for $l=0,\ldots,i$, 
for some $i\ge2$, we have
\begin{align*}
\y_{i+1}' &= \y_i' + \sum_{k=0}^ih_{i+1,k+1}(\x_{k+1}' - \y_k') \\
&= \y_i + \frac{2(i+1)}{i+2}(\x_{i+1} - \y_i) + \sum_{k=0}^{i-1}\paren{-\frac{2(k+1)}{(i+1)(i+2)}}(\x_{k+1} - \y_k) \\
&= \y_i + \paren{1+\frac{i}{i+2}}(\x_{i+1} - \y_i) 
	+ \frac{i}{i+2}\sum_{k=0}^{i-1}\paren{-\frac{2(k+1)}{i(i+1)}}(\x_{k+1} - \y_k) \\
&= \x_{i+1} + \frac{i}{i+2}(\x_{i+1} - \y_i) + \frac{i}{i+2}(\y_i + \y_{i-1} - 2\x_i) \\
&= \x_{i+1} + \frac{i}{i+2}(\x_{i+1} - \x_i) - \frac{i}{i+2}(\x_i - \y_{i-1}) = \y_{i+1},
\end{align*}
where the fourth equality uses 
\begin{align*}
\y_i &= \y_{i-1} + \frac{2i}{i+1}(\x_i - \y_{i-1})
+ \sum_{k=0}^{i-2}\left(-\frac{2(k+1)}{i(i+1)}\right)(\x_{k+1} - \y_k) \\
&= \y_{i-1} + 2(\x_i - \y_{i-1})
+ \sum_{k=0}^{i-1}\left(-\frac{2(k+1)}{i(i+1)}\right)(\x_{k+1} - \y_k)
.\end{align*}
\qed
\end{proof}
\end{proposition}

The proposed accelerated method has the inertia term $\frac{i}{i+2}(\x_{i+1} - \x_i)$,
similar to Nesterov's acceleration~\cite{nesterov:83:amf,nesterov:88:oaa} 
and G\"{u}ler's methods~\cite{guler:92:npp}.
However, the proposed method also has a correction term $- \frac{i}{i+2}(\x_i - \y_{i-1})$,
which is essential to guarantee an accelerated rate. 
Without such correction term, the accelerated method can diverge,
{\cred for} which we provide an example at the end of this section.
We leave further understanding the role of the proposed correction term as future work,
possibly via a differential equation perspective as in~\cite{su:16:ade}
for Nesterov's acceleration.
Note that a different correction term
for Nesterov's acceleration has been studied via the differential equation analysis for convex minimization
\cite{attouch:20:foo,shi:18:uta}.

The following theorem provides an accelerated rate of the proposed method
in terms of the fixed-point residual.\footnote{
The convergence of the fixed-point residual does not guarantee 
the convergence of the sequence of the iterates $\{\x_i\}$.
We leave analyzing the convergence of the sequence as future work,
possibly based on the convergence analysis in~\cite{chambolle:15:otc}
for Nesterov's fast gradient method~\cite{nesterov:83:amf,nesterov:88:oaa}
and FISTA~\cite{beck:09:afi} in convex minimization.
}

\begin{theorem}
\label{thm:appm,rate}
Let $\M\in\cM$ 
and let $\x_0,\y_0,\x_1,\y_1,\ldots\in\cH$ 
be generated by the proposed accelerated proximal point method.
Assume that $||\x_0 - \x_*||\le R$ for a constant $R>0$ 
and for some $\x_*\in X_*(\M)$. 
Then for any $i\ge1$,
\begin{align}
||\x_i - \y_{i-1}||^2 \le \frac{R^2}{i^2}
\label{eq:appm,rate}
.\end{align}
\begin{proof}
Using Lemma~\ref{lem:feas}, the general proximal point method with $\h$~\eqref{eq:hik} satisfies
\begin{align}
\frac{1}{R^2}||\x_N - \y_{N-1}||^2 \le \Bd_{D}(\h) \le \frac{1}{N^2}.
\label{eq:appm,bound}
\end{align}
Since the iterates of the method are recursive and do not depend on a given $N$,
the bound~\eqref{eq:appm,bound} generalizes to the intermediate iterates of the method.
By Proposition~\ref{prop:appm}, the proposed accelerated proximal point method also satisfies
the bound~\eqref{eq:appm,bound}, which concludes the proof.
\qed
\end{proof}
\end{theorem}

\begin{figure}[h!]
\begin{center}
\includegraphics[clip,width=0.48\textwidth]{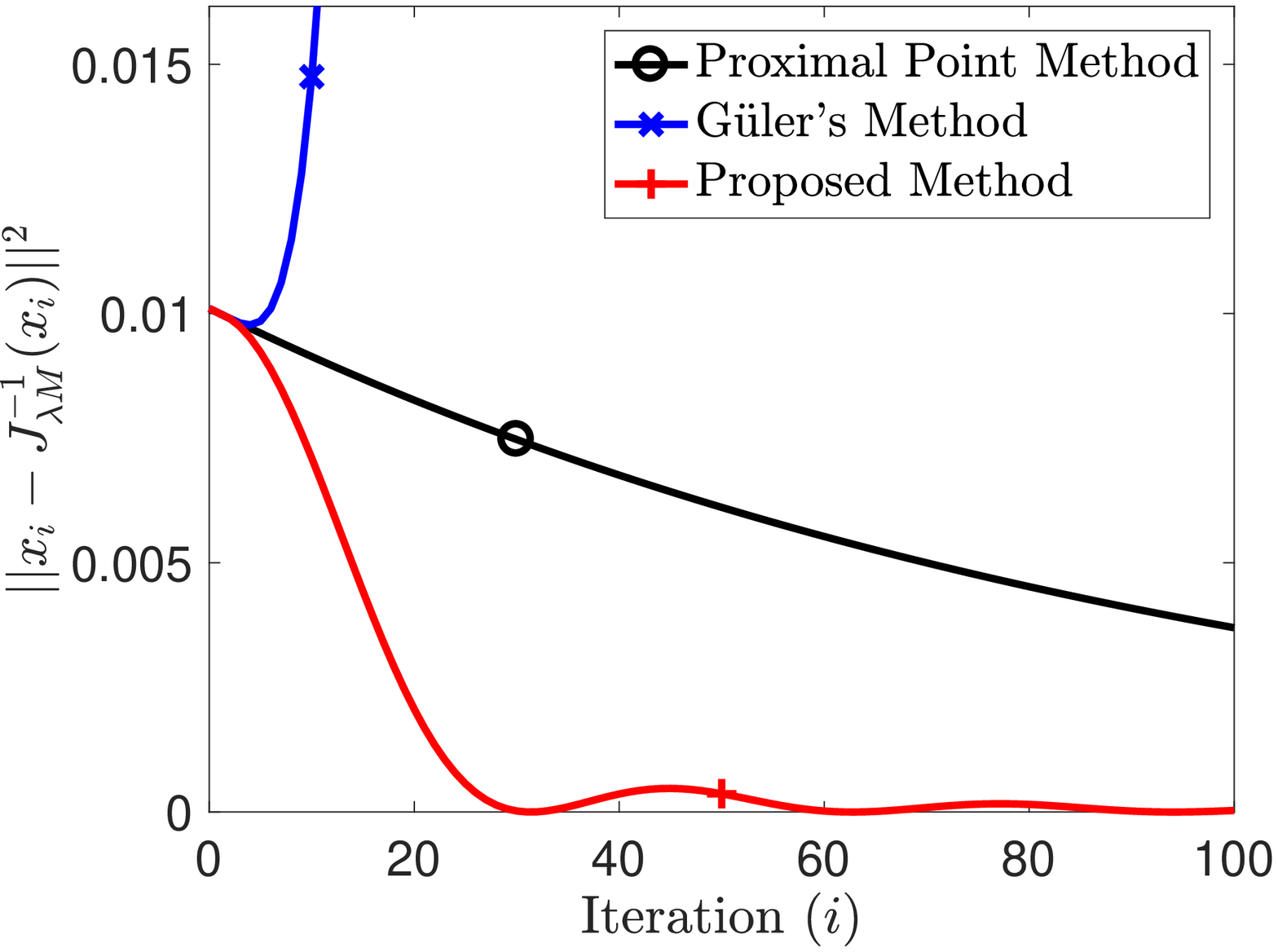}
\includegraphics[clip,width=0.48\textwidth]{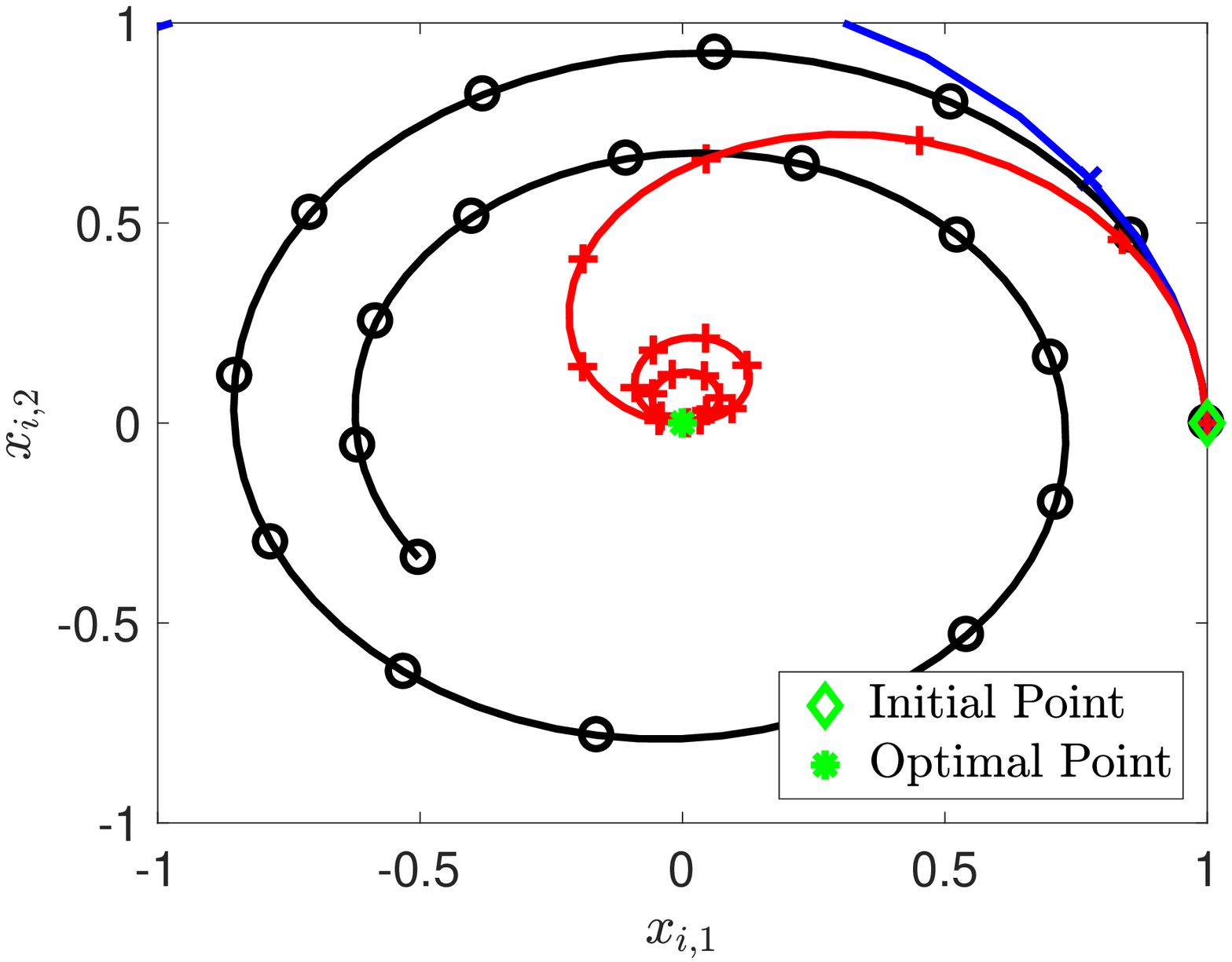}
\end{center}
\caption{Solving a worst-case monotone inclusion problem 
of the proximal point method with $\M$~\eqref{eq:M} with $N=100$;
(left) the fixed-point residual vs. iteration,
(right) the trajectory of the iterates $\x_i = [x_{i,1},\;x_{i,2}]^\top$ 
(markers are displayed every $5$th iterations).}
\label{fig:toy}
\end{figure}

The bound~\eqref{eq:ppm,rate} of the proximal point method 
was found to be exact in~\cite{gu:20:tsc} by specifying a certain operator $\M$ 
achieving the bound~\eqref{eq:ppm,rate} exactly; that is,
for given $N\ge2$, the proximal point method
exactly achieves the bound~\eqref{eq:ppm,rate}
for the operator
\begin{align}
\M\left[\begin{array}{c}
        u \\ v
        \end{array}\right]
= \frac{1}{\lambda\sqrt{N-1}}\left[\begin{array}{cc}
        0\; & \;1 \\
        -1\; & \;0
        \end{array}\right]
        \left[\begin{array}{c}
        u \\ v
        \end{array}\right]  
	\label{eq:M}
,\end{align}
with an initial point $\x_0 = [1\; 0]^\top$.
Such exact analysis is important since it reveals the worst-case behavior of the iterates of the method.
However, we were not able to show that the bound~\eqref{eq:appm,rate} 
of the proposed method is exact,
which we leave as future work.
Instead, we compared the behavior of the iterates of the proximal point method
and its accelerated variants on the operator $\M$ in~\eqref{eq:M}.
Figure~\ref{fig:toy} compares the proximal point method,
G\"{u}ler's first accelerated method with $\M$ instead of $\partial f$ (\ie, an instance of the inertia method)
and the proposed accelerated method,
with an initial point $\x_0 = [1\; 0]^\top$ and the optimal point $\x_* = \zero$.
Note that the G\"{u}ler's first method is almost equivalent to
the proposed accelerated method without the correction term $-\frac{i}{i+2}(\x_i - \y_{i-1})$,
and this exhibits diverging behavior in Fig.~\ref{fig:toy}.
The figure illustrates that
the correction term greatly helps the iterates to rapidly converge 
by reducing the radius of the orbit of the iterates,
compared to other methods.

We further investigate the behavior of the proposed method
for a convex-concave saddle-point problem
\begin{align}
\min_{\u\in\cH_1}\max_{\vv\in\cH_2}\; & \phi(\u,\vv)
\label{eq:minimax}
,\end{align}
where $\cH_1$ and $\cH_2$
denote real Hilbert spaces equipped with inner product $\inprod{\cdot}{\cdot}$,
and $\phi(\cdot,\vv)\in\cFone$, $-\phi(\u,\cdot)\in\cFtwo$,
which we further study in sections~\ref{sec:restart} and~\ref{sec:pmm}.
The saddle subdifferential of $\phi$,
\begin{align}
\left[\begin{array}{c}
\partial_{\u} \phi(\u,\vv) \\
\partial_{\vv} (-\phi(\u,\vv))
\end{array}\right]
\label{eq:saddle}
,\end{align}
is monotone~\cite{rockafellar:70:moa}.
The proposed accelerated method applied to~\eqref{eq:saddle}
with $\x_i := (\u_i,\vv_i)$ and $\x_* := (\u_*,\vv_*)$
(see Section~\ref{sec:pmm} for details)
%
satisfies
\begin{align}
\phi(\u_i,\vv_*) - \phi(\u_*,\vv_i) \le \frac{||\u_0 - \u_*||^2 + ||\vv_0 - \vv_*||^2}{4\lambda i}
\label{eq:phi}
\end{align}
for any $i\ge1$.
This is numerically conjectured by the PEP analysis in~\eqref{eq:Pfin}
with the objective function
$\frac{1}{R}||\x_N - \y_{N-1}||^2$
and the inequality
$\inprod{\x_N - \x_*}{\q_N} \ge 0$ in~\eqref{eq:Pfin}
replaced by $\phi(\u_N,\vv_*) - \phi(\u_*,\vv_N)$ 
and $\inprod{\x_N - \x_*}{\q_N} \ge \phi(\u_N,\vv_*) - \phi(\u_*,\vv_N)$, 
respectively.\footnote{
A convex-concave function $\phi$ satisfies
$\phi(\u_*,\vv_N) \ge \phi(\u_N,\vv_N) 
+ \inprod{\u_* - \u_N}{\q_{\u,N}}$
for $\q_{\u,N} \in \partial_{\u}\phi(\u_N,\vv_N)$
and $-\phi(\u_N,\vv_*) \ge -\phi(\u_N,\vv_N) 
+ \inprod{\vv_* - \vv_N}{-\q_{\vv,N}}$
for $-\q_{\vv,N} \in \partial_{\vv} (-\phi(\u_N,\vv_N))$.
Adding these two inequalities yields        
$\inprod{\x_N - \x_*}{\q_N} \ge \phi(\u_N,\vv_*) - \phi(\u_*,\vv_N)$,
where $\x_N := (\u_N,\vv_N)$ 
and $\q_N := (\q_{\u,N},-\q_{\vv,N})$.
}

\section{Restarting the accelerated proximal point method
for strongly monotone operators}
\label{sec:restart}

For strongly monotone operators,
the proximal point method has a linear rate
\eqref{eq:ppm,rate,strongly},
whereas the proposed accelerated method is not guaranteed to have such {\cred a} fast rate.
Technically, one should be able to find an accelerated method
for strong monotone operators via PEP,
as we did for the monotone operators in the previous section.
However, the resulting PEP problem, a reminiscent of~\eqref{eq:HD},
is much more difficult to solve, and we leave it as future work.
Instead, we consider a fixed restarting technique in 
\cite[Section 11.4]{nemirovski:94:emi}\cite[Section 5.1]{nesterov:13:gmf}
that restarts an accelerated method with a sublinear rate
every certain number of iterations
to yield a fast linear rate,
particularly for $\M\in\cMs$ in this section.

Suppose one restarts the proposed method every $k$ (inner) iterations
by initializing the $(j+1)$th outer iteration $\x_{j+1,0} = \y_{j+1,0} = \y_{j+1,-1}$
by $\x_{j,k}$,
where $\x_{j,l}$ and $\y_{j,l}$ denote iterates at the $j$th outer iteration
and $l$th inner iteration
for $j=0,1,\ldots$ and $l=-1,0,1,\ldots,k$.
Using the rate~\eqref{eq:appm,rate} (with $R = ||\x_{j,0} - \x_*||$) 
and the strong monotonicity condition~\eqref{eq:mono,strongly},
we have
\begin{align}
||\x_{j,k} - \y_{j,k-1}||^2 
	\le \frac{||\x_{j,0} - \x_*||^2}{k^2}
	\le \frac{1}{\mu^2 k^2}||\M\x_{j,0}||^2
\label{eq:M,restart}
\end{align}
for $j = 0,1,\ldots$.
Since $\frac{1}{\lambda}(\x_{j-1,k} - \y_{j-1,k-1}) \in \M\x_{j,0}$,
we have a linear rate
\begin{align}
||\x_{j,k} - \y_{j,k-1}||^2
        \le \frac{1}{\lambda^2\mu^2 k^2}||\x_{j-1,k} - \y_{j-1,k-1}||^2
\label{eq:restart}
.\end{align}
For a given $N=jk$ total number of steps,
{\cred
minimizing the overall rate with respect to $k$ yields
an optimal choice of the restarting interval given by}
$k_{\mathrm{opt}} {\cred\,\approx\,} \frac{e}{\lambda\mu}$,
where $e$ is Euler's number.
The corresponding linear rate is $O((e^{\lambda\mu/e})^{-2N})$.

We further investigate the behavior of the restarting {\cred technique}
for a saddle-point problem~\eqref{eq:minimax}
with an assumption that $\phi$ is strongly-convex-strongly-concave,
\ie, $\phi(\cdot,\vv) - \frac{\mu}{2}||\cdot||^2\in\cFone$ 
and $-\phi(\u,\cdot) - \frac{\mu}{2}||\cdot||^2\in\cFtwo$.
The associated saddle subdifferential~\eqref{eq:saddle} is $\mu$-strongly monotone.
For such case, using the rate~\eqref{eq:phi},
and the inequalities $\phi(\u_*,\vv_*) + \frac{\mu}{2}||\u - \u_*||^2 \le \phi(\u,\vv_*)$
and $-\phi(\u_*,\vv_*) + \frac{\mu}{2}||\vv - \vv_*||^2 \le -\phi(\u_*,\vv)$,
the proposed method with restarting every $k$ iterations satisfies
\begin{align}
\phi(\u_{j,k},\vv_*) - \phi(\u_*,\vv_{j,k}) 
\le \frac{1}{2\lambda\mu k}(\phi(\u_{j,0},\vv_*) - \phi(\u_*,\vv_{j,0}))
\label{eq:phi,restart}
\end{align}
for $j=0,1,\ldots$.
The associated optimal restarting interval is 
$k_{\mathrm{opt}}^\phi {\cred\,\approx\,} \frac{e}{2\lambda\mu}$,
which is twice smaller than $k_{\mathrm{opt}}$.
The corresponding linear rate is also $O((e^{\lambda\mu/e})^{-2N})$,
whereas the proximal point method has the rate $O((1+\lambda\mu)^{-2N})$
in~\eqref{eq:ppm,rate,strongly}.
{\cred
For any given positive $\mu$,
there is no positive $\lambda$
that satisfies both $k_{\mathrm{opt}}^\phi(\lambda,\mu)\ge1$ and $e^{\lambda\mu/e} > 1+\lambda\mu$.
This contrasts with the fact 
that the worst-case rate of optimally restarting the proposed method 
is not slower than that of the proximal point method.
This implies that the bounds~\eqref{eq:restart} and~\eqref{eq:phi,restart} are not exact,
and we leave finding their tight bounds
%
as future work.
}
The numerical experiment below (and those in Section~\ref{sec:appl})
suggests that the restarting technique can perform better 
than the proximal point method.

We consider a toy problem that is a combination of the worst-case problems
in $\cM$ and $\cMs$ for the proximal point method:
\begin{align}
\M\left[\begin{array}{c}
        u \\ v
        \end{array}\right]
= \left(\frac{1}{\lambda\sqrt{N-1}}\left[\begin{array}{cc}
        0\; & \;1 \\
        -1\; & \;0
        \end{array}\right]
	+ \left[\begin{array}{cc}
        \mu\; & \;0 \\
        0\; & \;\mu
        \end{array}\right]
	\right)
        \left[\begin{array}{c}
        u \\ v
        \end{array}\right]  
        \label{eq:Mb}
,\end{align}
which is the saddle subdifferential operator of 
$\phi(u,v) = \frac{\mu}{2}u^2 + \frac{1}{\lambda\sqrt{N-1}}uv - \frac{\mu}{2}v^2$.
We choose
$N = 100$, $\lambda = 1$ and $\mu = 0.02$.
The optimal restarting intervals are
$k_{\mathrm{opt}} \approx 136$
and $k_{\mathrm{opt}}^\phi \approx 68$,
and we run $200$ iterations in the experiment,
where restarting intervals $17$, $34$, $68$, and $136$ are considered.
Figure~\ref{fig:toy,strongly} compares the proximal point method, its accelerated variants,
and the proposed accelerated method with restarting,
with an initial point $\x_0 = [1\; 0]^\top$ and the optimal point $\x_* = \zero$.
Figure~\ref{fig:toy,strongly} presents that
the proximal point method has a linear rate that is faster than the proposed method
(with a sublinear rate),
while the restarting greatly accelerates the proposed method with a fast linear rate.
Figure~\ref{fig:toy,strongly} {\cred also} illustrates 
that the optimal restarting intervals $k_{\mathrm{opt}}$ and $k_{\mathrm{opt}}^\phi$
{\cred for strongly monotone operators and strongly-convex-strongly-concave functions, respectively,
are not optimal for this specific case.}
Examples in the next section also present
that the restarting can be useful even without strong monotonicity
(but possibly with local strong monotonicity).

\begin{figure}[h!]
\begin{center}
\includegraphics[clip,width=0.48\textwidth]{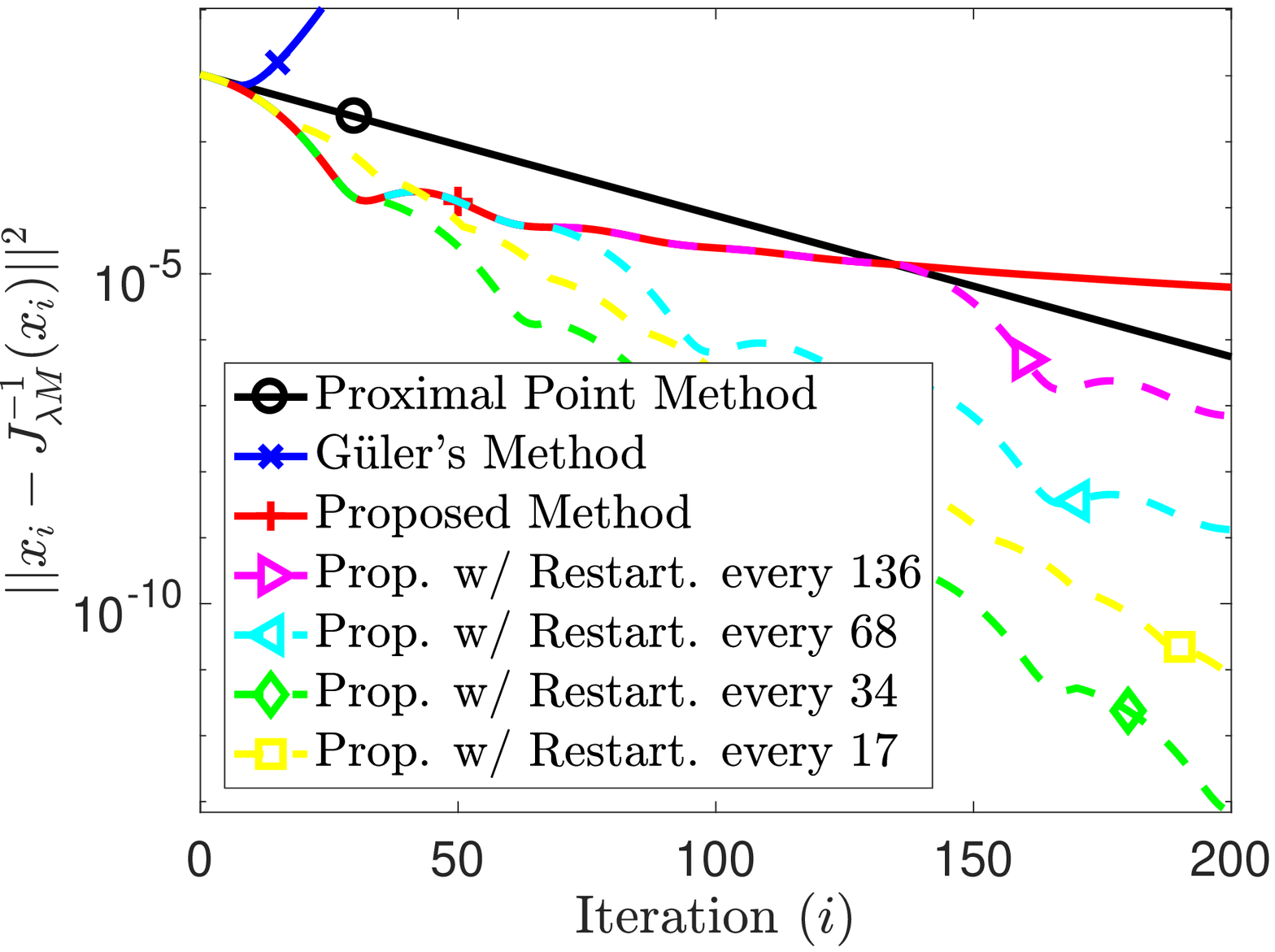}
\includegraphics[clip,width=0.48\textwidth]{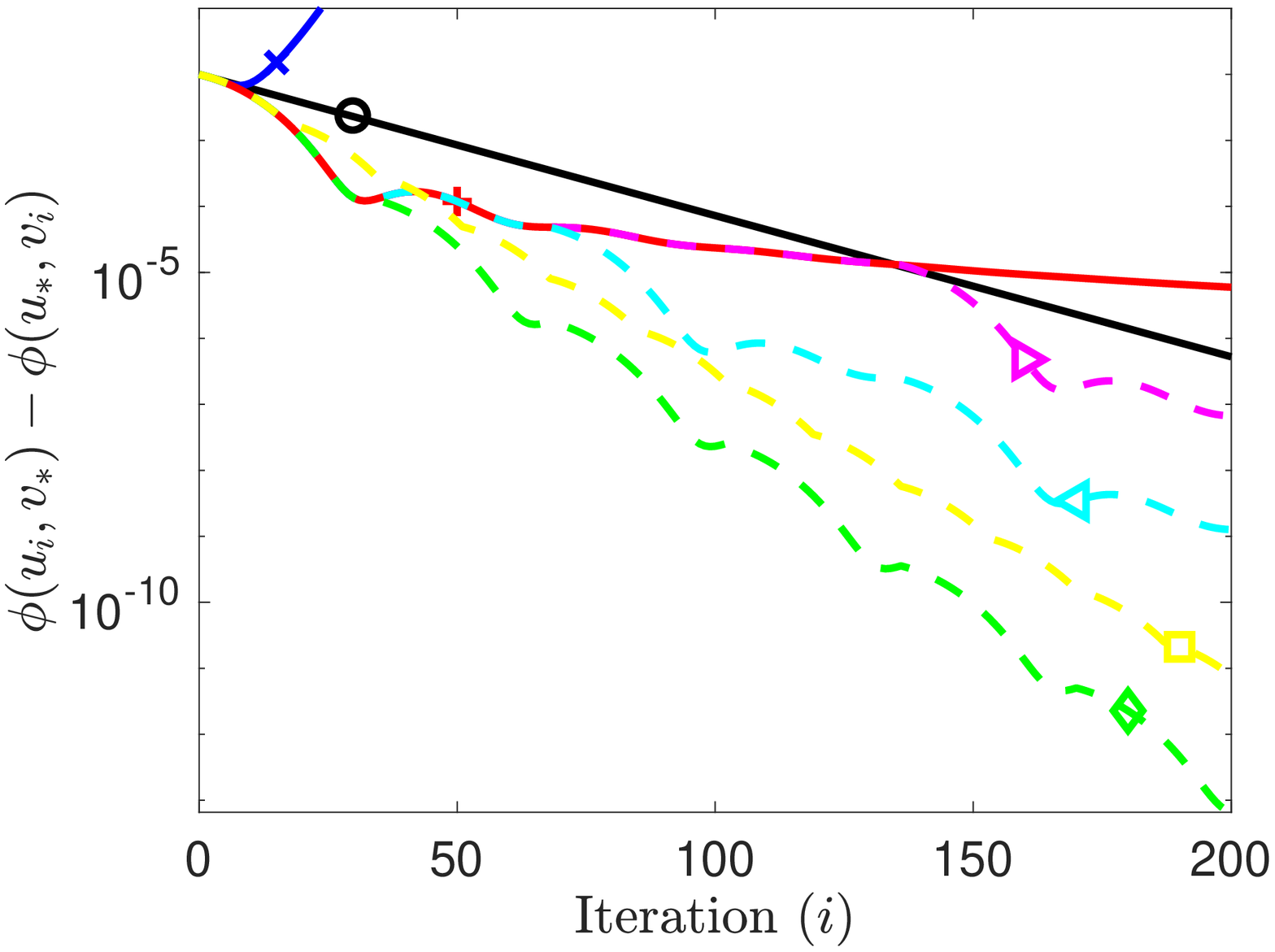}
\end{center}
\caption{Solving a strongly monotone inclusion problem
with $\M$~\eqref{eq:Mb};
(left) the fixed-point residual vs. iteration,
(right) the function residual vs. iteration.}
\label{fig:toy,strongly}
\end{figure}

\section{Applications of the accelerated proximal point method}
\label{sec:appl}

As mentioned earlier,
the proximal point method for maximally monotone operators
include various well-known convex optimization methods.
These include
the augmented Lagrangian (\ie, the method of multipliers),
the proximal method of multipliers,
and ADMM.
The augmented Lagrangian method
is equivalent to the proximal point method 
directly solving the dual convex minimization problem~\cite{rockafellar:76:ala},
so G\"{u}ler's methods~\cite{guler:92:npp} already provide acceleration,
whereas other instances of the proximal point method have no known accelerations yet.
Thus, this section introduces accelerations
to well-known instances of the proximal point method,
which were not possible previously to the best of our knowledge
(under this paper's setting).

\subsection{Accelerating the proximal point method for convex-concave saddle-point problem}
\label{sec:pmm}


This section considers a convex-concave saddle-point problem~\eqref{eq:minimax},
where the associated saddle subdifferential operator~\eqref{eq:saddle}
is monotone.
\cite{rockafellar:76:moa} applied the proximal point method on such operator
to solve the convex-concave saddle-point problem,
and this section further applies the proposed acceleration to such proximal point method as below.

\fbox{
\begin{minipage}[t]{0.85\linewidth}
\vspace{-10pt}
\begin{flalign*}
&\quad \text{\bf Accelerated Proximal Point Method for Convex-Concave Saddle-Point Problem} & \\
&\qquad \text{Input: } \phi(\cdot,\vv)\in\cFone,\; -\phi(\u,\cdot)\in\cFtwo,\; 
	\hat{\u}_0\in\cH_1,\; \hat{\vv}_0\in\cH_2, & \\ 
&\qquad \qquad\quad	
	\x_0 = \y_0 = \y_{-1} = (\hat{\u}_0, \hat{\vv}_0),\; \lambda \in\reals_{++}. & \\
&\qquad \text{For } i = 0,1,\ldots & \\
&\qquad \qquad \x_{i+1} = (\u_{i+1},\vv_{i+1}) = \arg\min_{\u\in\cH_1}\max_{\vv\in\cH_2} \left\{\phi(\u,\vv) 
                + \frac{1}{2\lambda}||\u - \hat{\u}_i||^2 - \frac{1}{2\lambda}||\vv - \hat{\vv}_i||^2\right\}, & \\
&\qquad \qquad \y_{i+1} = (\hat{\u}_{i+1},\hat{\vv}_{i+1}) = \x_{i+1} + \frac{i}{i+2}(\x_{i+1} - \x_i) - \frac{i}{i+2}(\x_i - \y_{i-1}). &
\end{flalign*}
\end{minipage}
} \vspace{5pt}

One primary use of this accelerated method
is the following convex-concave Lagrangian problem
\begin{align}
\min_{\u\in\cH_1}\max_{\vv\in\cH_2}\; & \left\{L(\u,\vv) := f(\u) + \inprod{\vv}{\A\u - \bb}\right\},
\label{eq:minmaxL}
\end{align}
associated with the linearly constrained problem
\begin{align}
\min_{\u\in\cH_1}\; & f(\u) \\
\text{subject to }\; & \A\u = \bb, \nonumber
\end{align}
where $\A\in\cB(\cH_1,\cH_2)$ 
and $\bb\in\cH_2$.
The resulting method
is called the proximal method of multipliers in~\cite{rockafellar:76:ala},
and applying the proposed acceleration to this method leads to below.

\fbox{
\begin{minipage}[t]{0.85\linewidth}
\vspace{-10pt}
\begin{flalign*}
&\quad \text{\bf Accelerated Proximal Method of Multipliers} & \\
&\qquad \text{Input: } \phi(\cdot,\vv)\in\cFone,\; -\phi(\u,\cdot)\in\cFtwo,\; 
        \hat{\u}_0\in\cH_1,\; \hat{\vv}_0\in\cH_2, & \\
&\qquad\qquad\quad 
        \x_0 = \y_0 = \y_{-1} = (\hat{\u}_0, \hat{\vv}_0),\; \lambda \in\reals_{++}. & \\
&\qquad \text{For } i = 0,1,\ldots & \\
&\qquad \qquad \u_{i+1} = \arg\min_{\u\in\cH_1} \left\{L(\u,\hat{\vv}_i) + \frac{\lambda}{2}||\A\u - \bb||^2
                + \frac{1}{2\lambda}||\u - \hat{\u}_i||^2\right\}, & \\
&\qquad \qquad \x_{i+1} = (\u_{i+1},\hat{\vv}_i + \lambda(\A\u_{i+1} - \bb)), & \\
&\qquad \qquad \y_{i+1} = (\hat{\u}_{i+1},\hat{\vv}_{i+1}) = \x_{i+1} + \frac{i}{i+2}(\x_{i+1} - \x_i) - \frac{i}{i+2}(\x_i - \y_{i-1}). &
\end{flalign*}
\end{minipage}
} \vspace{5pt}

Note that this method without the acceleration and the term $\frac{1}{2\lambda}||\u - \u_i||^2$ 
reduces to the augmented Lagrangian method.
This method has an advantage over the augmented Lagrangian method and its accelerated variants;
the primal iterate $\u_{i+1}$ is uniquely defined 
with a better conditioning.

\begin{example}
\label{ex:basis_pursuit}
We apply the accelerated proximal method of multipliers
to a basis pursuit problem
\begin{align}
\min_{\u\in\reals^{d_1}}\; & ||\u||_1 \label{eq:basis_pursuit} \\
\text{subject to }\; & \A\u = \bb, \nonumber
\end{align}
where $\A\in\reals^{d_2\times d_1}$ and $\bb\in\reals^{d_2}$.
In the experiment,
we choose $d_1 = 100$, $d_2 = 20$,
and randomly generated $\A$.
A true sparse $\u_{\mathrm{true}}$ is randomly generated followed by a thresholding to sparsify nonzero elements,
and $\bb$ is then given by $\A\u_{\mathrm{true}}$.
We run $100$ iterations of the proximal method of multipliers and its variants
with $\lambda = 0.01$ and initial $\x_0 = \zero$.
Since the $\u_{i+1}$-update does not have a closed form,
we used a sufficient number of iterations 
to solve the $\u_{i+1}$-update
using the strongly convex version of FISTA~\cite{beck:09:afi}
in~\cite[Theorem 4.10]{chambolle:16:ait}.

Figure~\ref{fig:basis} compares the proximal method of multipliers
and its accelerated variants.
Similar to Fig.~\ref{fig:toy}, G\"{u}ler's first accelerated version diverges,
while the proposed method has accelerating behavior, compared to the non-accelerated version.
The proposed method exhibits an oscillation in Fig.~\ref{fig:basis}
(and a subtle oscillation in Fig.~\ref{fig:toy}).
This might be due to \emph{high momentum}, owing from the acceleration,
discussed in~\cite{odonoghue:15:arf}.
So in Fig.~\ref{fig:basis}
we heuristically restarted the method every $30$ iterations
to avoid such oscillation and accelerate,
as suggested in~\cite{odonoghue:15:arf}.
Developing an approach to appropriately choosing a restarting interval
or adaptively restarting the method as in~\cite{odonoghue:15:arf}
for such problem
are left as future work.\footnote{
We found that adaptively restarting the method when the fixed-point residual increases
seems to be a good option in practice.
}

\begin{figure}[h!]
\begin{center}
\includegraphics[clip,width=0.55\textwidth]{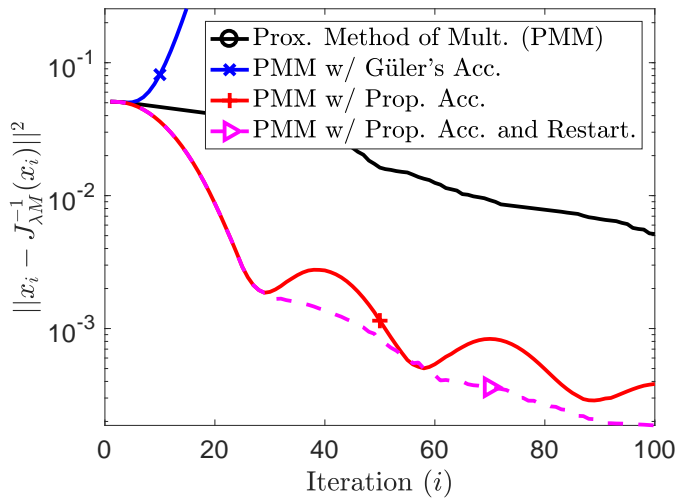}
\end{center}
\caption{Solving a basis pursuit problem~\eqref{eq:basis_pursuit};
the fixed-point residual vs. iteration.}
\label{fig:basis}
\end{figure}

\end{example}

\subsection{Accelerating the primal-dual hybrid gradient method}
\label{sec:acpm}

This section considers a linearly coupled convex-concave saddle-point problem
\begin{align}
\min_{\u\in\cH_1}\max_{\vv\in\cH_2} \left\{\phi(\u,\vv) \equiv f(\u) + \inprod{\K\u}{\vv} - g(\vv)\right\}
,\end{align}
where $f\in\cFone$, $g\in\cFtwo$ and $\K\in\cB(\cH_1,\cH_2)$.
One widely known method for such problem
is the primal-dual hybrid gradient (PDHG) method~\cite{chambolle:11:afo,esser:10:agf},
which is a \emph{preconditioned} proximal point method (with $\lambda = 1$)
for the saddle subdifferential operator of $\phi$~\eqref{eq:saddle}
\cite{chambolle:16:ote,he:12:cao}.
The associated preconditioner is
\begin{align}
\P = \left[\begin{array}{cc}
        \frac{1}{\tau}\I \;&\; -\K^* \\
        -\K \;&\; \frac{1}{\sigma}\I
        \end{array}\right]
,\end{align} 
which is positive definite when $\tau\sigma||\K||^2<1$,
where $||\K|| = \sup_{||\x||\le1} ||\K\x||$.
As mentioned in remark~\ref{remark:precond},
we can directly apply our results to the PDHG method as below.

\fbox{
\begin{minipage}[t]{0.85\linewidth}
\vspace{-10pt}
\begin{flalign*} 
&\quad \text{\bf Accelerated PDHG Method} & \\
&\qquad \text{Input: } f\in\cFone,\; g\in\cFtwo,\; \K\in\cB(\cH_1,\cH_2),\;
                \hat{\u}_0\in\cH_1,\; \hat{\vv}_0\in\cH_2,\;
                \tau\sigma||\K||^2<1, & \\
&\qquad \qquad\quad \x_0 = \y_0 = \y_{-1} = (\hat{\u}_0,\hat{\vv}_0). & \\
&\qquad \text{For } i = 0,1,\ldots & \\
&\qquad \qquad \u_{i+1} = \argmin{\u\in\cH_1} \left\{f(\u) + \inprod{\K\u}{\hat{\vv}_i} + \frac{1}{2\tau}||\u - \hat{\u}_i||^2\right\} & \\
&\qquad \qquad \vv_{i+1} = \argmin{\vv\in\cH_2} \left\{g(\vv) - \inprod{\K(2\u_{i+1} - \hat{\u}_i)}{\vv}
                        + \frac{1}{2\sigma}||\vv - \hat{\vv}_i||^2\right\} & \\
&\qquad \qquad \x_{i+1} = (\u_{i+1},\vv_{i+1}) & \\
&\qquad \qquad \y_{i+1} = (\hat{\u}_{i+1},\hat{\vv}_{i+1}) = \x_{i+1} + \frac{i}{i+2}(\x_{i+1} - \x_i)
                        - \frac{i}{i+2}(\x_i - \y_{i-1}) &
\end{flalign*}
\end{minipage}
} \vspace{5pt}

\begin{corollary}
Assume that $\inprod{\P(\x_0 - \x_*)}{\x_0 - \x_*} \le R^2$ 
for some $\x_* \in X_*(\M)$.
The PDHG method satisfies
\begin{align*}
\inprod{\P(\x_i - \x_{i-1})}{\x_i - \x_{i-1}} \le \left(1 - \frac{1}{i}\right)^{i-1}\frac{R^2}{i}
,\end{align*}
and the proposed accelerated PDHG method satisfies
\begin{align*}
\inprod{\P(\x_i - \y_{i-1})}{\x_i - \y_{i-1}} \le \frac{R^2}{i^2}
.\end{align*}
\end{corollary}

\begin{example}
We apply the accelerated PDHG method 
to the bilinear game problem
\begin{align}
\min_{\u\in\reals^{d_1}}\max_{\vv\in\reals^{d_2}} \inprod{\aa}{\u} + \inprod{\K\u}{\vv} - \inprod{\bb}{\vv}
\label{eq:matrix_game}
,\end{align}
where $\K\in\reals^{d_2\times d_1}$, $\aa\in\reals^{d_1}$ and $\bb\in\reals^{d_2}$.
The main part of the corresponding method is as below:
\begin{align}
\u_{i+1} &= \hat{\u}_i - \tau(\K^*\hat{\vv}_i + \aa) \\
\vv_{i+1} &= \hat{\vv}_i + \sigma(\K(2\u_{i+1} - \hat{\u}_i) - \bb) \nonumber
.\end{align}

In the experiment, we choose $d_1 = 1000$, $d_2 = 500$, and a matrix $\K$ and vectors $\aa,\bb$ are randomly generated.
We run $100$ iterations of the PDHG method and its variants 
with initial $\hat{\u}_0 = [10\;\cdots\;10]^\top$, $\hat{\vv}_0 = [10\;\cdots\;10]^\top$ and $\tau=\sigma=\frac{0.99}{||\K||}$.
Figure~\ref{fig:cpm} plots the \emph{preconditioned} fixed-point residual,
where G\"{u}ler's first accelerated method diverges.
The PDHG method and its proposed accelerated variant are comparable in this experiment,
and heuristically restarting the accelerated method every $10$ iterations yields a big acceleration.
While~\cite{chambolle:11:afo} found restarting (reinitializing)
a relaxed PDHG method not useful,
our experiment suggests that
restarting can be effective in some practical cases.

\begin{figure}[h!]
\begin{center}
\includegraphics[clip,width=0.55\textwidth]{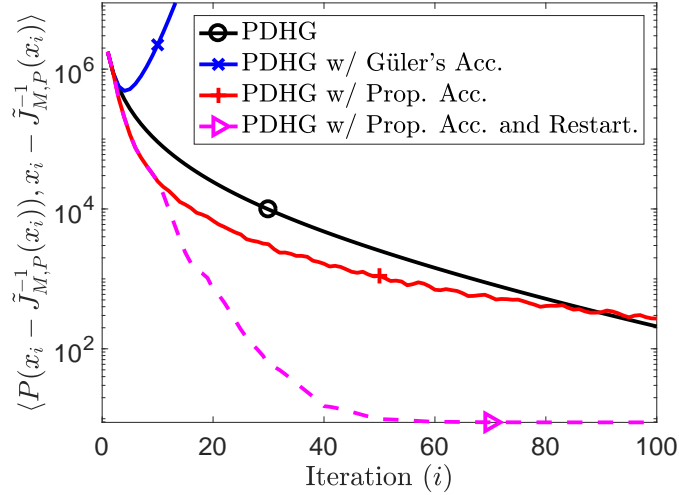}
\end{center}
\caption{
Solving a bilinear game problem~\eqref{eq:matrix_game}
by a \emph{preconditioned} proximal point method;
the \emph{preconditioned} fixed-point residual vs. iteration.
$\tilde{\J}_{\M,\P} := (\P + \M)^{-1}\P$ denotes the \emph{preconditioned} resolvent operator.}
\label{fig:cpm}
\end{figure}

\end{example}

\subsection{Accelerating the Douglas-Rachford splitting method}

This section considers a monotone inclusion problem
in a form
\begin{align}
\Find \;\; \x\in\cH \quad \text{subject to} \quad \zero \in (\M_1 + \M_2)\x 
\label{eq:mono,split}
\end{align}
for $\M_1,\M_2\in\cM$,
where $\J_{\rho\M_1}$ and $\J_{\rho\M_2}$ are more efficient
than $\J_{\rho(\M_1 + \M_2)}$
for a positive real number $\rho$.
For such problem,
the Douglas-Rachford splitting method~\cite{douglas:56:otn,lions:79:saf}
that iteratively applies the operator
\begin{align}
\G_{\rho,\M_1,\M_2} := \J_{\rho\M_1}\circ(2\J_{\rho\M_2} - \I) + (\I - \J_{\rho\M_2})
\label{eq:dg}
\end{align}
has been found to be effective
in many applications
including ADMM, which we discuss in the next section.

In~\cite[Theorem 4]{eckstein:92:otd},
the Douglas-Rachford operator~\eqref{eq:dg}
was found to be a resolvent $\J_{\M_{\rho,\M_1,\M_2}}$ of a maximally monotone operator
\begin{align}
\M_{\rho,\M_1,\M_2} := \G_{\rho,\M_1,\M_2}^{-1} - \I
\label{eq:Mrho}
.\end{align}
In other words, the Douglas-Rachford splitting method is an instance 
of the proximal point method (with $\lambda = 1$) as
\begin{align}
\nuu_{i+1} = \J_{\M_{\rho,\M_1,\M_2}}(\nuu_i) = \G_{\rho,\M_1,\M_2}(\nuu_i)
\end{align}
for $i=0,1,\ldots$.
Therefore, we can apply the proposed acceleration
to the Douglas-Rachford splitting method as below.

\fbox{
\begin{minipage}[t]{0.85\linewidth}
\vspace{-10pt}
\begin{flalign*} 
&\quad \text{\bf Accelerated Douglas-Rachford Splitting Method} & \\
&\qquad \text{Input: } \M_1,\M_2\in\cM,\; \nuu_0 = \etaa_0 = \etaa_{-1}\in\cH,\; \rho\in\reals_{++}. & \\
&\qquad \text{For } i = 0,1,\ldots & \\
&\qquad \qquad \nuu_{i+1} = \G_{\rho,\M_1,\M_2}(\etaa_i) & \\
&\qquad \qquad \etaa_{i+1} = \nuu_{i+1} + \frac{i}{i+2}(\nuu_{i+1} - \nuu_i) 
                        - \frac{i}{i+2}(\nuu_i - \etaa_{i-1}) &
\end{flalign*}
\end{minipage}
} \vspace{5pt}

Using~\eqref{eq:ppm,rate} and~\eqref{eq:appm,rate},
we have the following worst-case rates
for the Douglas-Rachford splitting method and its accelerated variant.
Finding exact bounds for the Douglas-Rachford splitting method and its variant
is left as future work;
\cite{ryu:20:osp} used PEP to analyze the exact worst-case rate of Douglas-Rachford splitting method
under some additional conditions.

\begin{corollary}
Assume that $||\nuu_0 - \nuu_*||\le R$ for some $\nuu_*\in X_*(\M_{\rho,\M_1,\M_2})$.
The Douglas-Rachford splitting method satisfies
\begin{align}
||\nuu_i - \etaa_{i-1}||^2 \le \left(1 - \frac{1}{i}\right)^{i-1}\frac{R^2}{i}
,\end{align}
and the proposed accelerated Douglas-Rachford splitting method satisfies
\begin{align}
||\nuu_i - \etaa_{i-1}||^2 \le \frac{R^2}{i^2}
.\end{align}
\end{corollary}

\cite{eckstein:88:tlm,eckstein:92:otd}
illustrated that ADMM is equivalent to the Douglas-Rachford splitting method on the dual problem,
so we naturally develop an accelerated ADMM in the next section
and provide numerical experiment of the accelerated ADMM 
and thus the accelerated Douglas-Rachford splitting method.

\subsection{Accelerating the alternating direction method of multipliers (ADMM)}

Let $\cH_1,\cH_2,\cG$ be real Hilbert spaces equipped with inner product $\inprod{\cdot}{\cdot}$.
This section considers a linearly constrained convex problem
\begin{align}
\min_{\x\in\cH_1,\z\in\cH_2}\; & f(\x) + g(\z) \label{eq:primal} \\
\text{subject to }\; & \A\x + \B\z = \cc, \nonumber
\end{align}
where $f\in\cFone$, $g\in\cFtwo$, $\A\in\cB(\cH_1,\cG)$, $\B\in\cB(\cH_2,\cG)$ and $\cc \in\cG$.
Its dual problem is
\begin{align}
\max_{\nuu\in\cG}\; & \left\{- f^*(-\A^*\nuu) - g^*(-\B^*\nuu) + \inprod{\cc}{\nuu}\right\}
\label{eq:dual}
,\end{align}
where $f^*(\y) := \sup_{\x\in\cH_1} \{\inprod{\y}{\x} - f(\x)\}$ 
and $g^*(\y) := \sup_{\z\in\cH_2} \{\inprod{\y}{\z} - g(\z)\}$
are the conjugate functions of $f$ and $g$, respectively.
The dual problem~\eqref{eq:dual} is equivalent to the following monotone inclusion problem
\begin{align}
\Find \;\; \nuu\in\cG \quad \text{subject to} \quad 
	\zero \in -\A\partial f^*(-\A^*\nuu) - \B\partial g^*(-\B^*\nuu) - \cc
\label{eq:dual,mono}
.\end{align}
We next use the connection between ADMM for solving~\eqref{eq:primal}
and the Douglas-Rachford splitting method for solving~\eqref{eq:dual,mono}
in~\cite[Proposition 9]{davis:16:cra}\cite{ryu:16:apo}
to develop an accelerated ADMM,
using the accelerated Douglas-Rachford splitting method in the previous section.

Denoting
\begin{align}
\M_1 &:= -\A\partial f^*(-\A^*\cdot) - \cc
\quad \text{and} \quad
\M_2 := -\B\partial g^*(-\B^*\cdot)
\label{eq:M1M2}
\end{align}
converts the problem~\eqref{eq:dual,mono} into a form of the monotone inclusion problem~\eqref{eq:mono,split}.
Then we use the following equivalent form of 
the accelerated Douglas-Rachford splitting method
to solve~\eqref{eq:mono,split} with~\eqref{eq:M1M2}:
\begin{align}
\zet_{i+1} &= \J_{\rho\M_2}(\etaa_i) \\
\xii_{i+1} &= \J_{\rho\M_1}(2\zet_{i+1} - \etaa_i) \nonumber \\
\nuu_{i+1} &= \etaa_i + (\xii_{i+1} - \zet_i) \nonumber \\
\etaa_{i+1} & = \nuu_{i+1} + \frac{i}{i+2}(\nuu_{i+1} - \nuu_i) - \frac{i}{i+2}(\nuu_i - \etaa_{i-1}) \nonumber
\end{align}
for $i=0,1,\ldots$.
Replacing the resolvent operators of $\M_1$ and $\M_2$ in~\eqref{eq:M1M2} by minimization steps yields
\begin{align}
\z_{i+1} &= \argmin{\z\in\cH_2} \left\{g(\z) + \inprod{\etaa_i}{\B\z} 
		+ \frac{\rho}{2}||\B\z||^2\right\} \\
\zet_{i+1} &= \etaa_i + \rho\B\z_{i+1} \nonumber \\
\tilde{\x}_{i+1} &= \argmin{\x\in\cH_1} \left\{f(\x) + \inprod{\etaa_i + 2\rho\B\z_{i+1}}{\A\x - \cc}
		+ \frac{\rho}{2}||\A\x - \cc||^2\right\} \nonumber \\
\xii_{i+1} &= \etaa_i + \rho(\A\tilde{\x}_{i+1} - \cc) + 2\rho\B\z_{i+1} \nonumber \\
\nuu_{i+1} &= \etaa_i + \rho(\A\tilde{\x}_{i+1} + \B\z_{i+1} - \cc) \nonumber \\
\etaa_{i+1} & = \nuu_{i+1} + \frac{i}{i+2}(\nuu_{i+1} - \nuu_i) 
		- \frac{i}{i+2}(\nuu_i - \etaa_{i-1}). \nonumber
\end{align}
By discarding $\zet_i$ and $\xii_i$, and defining
\begin{align}
\hat{\nuu}_i := \nuu_i - \rho(\A\tilde{\x}_i - \cc)
\quad\text{and}\quad
\hat{\etaa}_i := \etaa_i - \rho(\A\tilde{\x}_i - \cc), 
\end{align}
for $i=0,1,\ldots$,
we have
\begin{align}
\z_{i+1} &= \argmin{\z\in\cH_2} \left\{g(\z) + \inprod{\hat{\etaa}_i + \rho(\A\tilde{\x}_i - \cc)}{\B\z} 
                + \frac{\rho}{2}||\B\z||^2\right\} \\
		&= \argmin{\z\in\cH_2} \left\{g(\z) + \inprod{\hat{\etaa}_i}{\A\tilde{\x}_i + \B\z - \cc} 
                + \frac{\rho}{2}||\A\tilde{\x}_i + \B\z- \cc||^2\right\} \nonumber \\
\tilde{\x}_{i+1} &= \argmin{\x\in\cH_1} \left\{f(\x) 
                + \inprod{\hat{\nuu}_{i+1} + \rho\B\z_{i+1}}{\A\x - \cc}
		+ \frac{\rho}{2}||\A\x - \cc||^2\right\} \nonumber \\
		&= \argmin{\x\in\cH_1} \left\{f(\x) 
                + \inprod{\hat{\nuu}_{i+1}}{\A\x + \B\z_{i+1} - \cc}
                + \frac{\rho}{2}||\A\x + \B\z_{i+1} - \cc||^2\right\} \nonumber \\
\hat{\nuu}_{i+1} &= \hat{\etaa}_i + \rho(\A\tilde{\x}_i + \B\z_{i+1} - \cc) \nonumber \\
\hat{\etaa}_{i+1} & = \hat{\nuu}_{i+1} + \frac{i}{i+2}(\hat{\nuu}_{i+1} - \hat{\nuu}_i + \rho\A(\tilde{\x}_{i+1} - \tilde{\x}_i)) 
                - \frac{i}{i+2}(\hat{\nuu}_i - \hat{\etaa}_{i-1} + \rho\A(\tilde{\x}_i - \tilde{\x}_{i-1})) \nonumber
.\end{align}
Then, replacing $\tilde{\x}_i$ by $\x_{i+1}$ and reordering steps appropriately
yield the following accelerated version of ADMM,
which reduces to the standard ADMM 
when we let $\hat{\etaa}_i = \hat{\nuu}_i$ for $i=0,1,\ldots$.

\fbox{
\begin{minipage}[t]{0.85\linewidth}
\vspace{-10pt}
\begin{flalign*}
&\quad \text{\bf Accelerated Alternating Direction Method of Multipliers} & \\
&\qquad \text{Input: } f\in\cFone,\; g\in\cFtwo,\; \A\in\cB(\cH_1,\cG),\; \B\in\cB(\cH_2,\cG),\;
		\x_0\in\cH_1,\; \z_0\in\cH_2,\;
		\hat{\nuu}_0\in\cG,\; \rho\in\reals_{++}. & \\
&\qquad \text{For } i = 0,1,\ldots & \\
&\qquad \qquad \x_{i+1} = \argmin{\x\in\cH_1} \left\{f(\x) + \inprod{\hat{\nuu}_i}{\A\x +\B\z_i - \cc}
                + \frac{\rho}{2}||\A\x + \B\z_i - \cc||^2\right\} & \\
&\qquad \qquad \hat{\etaa}_i = \begin{cases}
		\hat{\nuu}_i & i = 0,1, \\
		\hat{\nuu}_i + \frac{i-1}{i+1}(\hat{\nuu}_i - \hat{\nuu}_{i-1} + \rho\A(\x_{i+1} - \x_i)) 
                - \frac{i-1}{i+1}(\hat{\nuu}_{i-1} - \hat{\etaa}_{i-2} + \rho\A(\x_i - \x_{i-1})),
		& i = 2,3,\ldots
		\end{cases}
		& \\
&\qquad \qquad \z_{i+1} = \argmin{\z\in\cH_2} \left\{g(\z) + \inprod{\hat{\etaa}_i}{\A\x_{i+1} + \B\z - \cc} 
                + \frac{\rho}{2}||\A\x_{i+1} + \B\z - \cc||^2\right\} & \\
&\qquad \qquad \hat{\nuu}_{i+1} = \hat{\etaa}_i + \rho(\A\x_{i+1} + \B\z_{i+1} - \cc) &
\end{flalign*}
\end{minipage}
} \vspace{5pt}

Since 
\begin{align}
\nuu_i - \etaa_{i-1}
= \hat{\nuu}_i - \hat{\etaa}_{i-1} - \rho(\A\x_{i+1} - \A\x_i) 
= \rho(\A\x_{i+1} + \B\z_i - \cc), 
\end{align}
we have the following worst-case rates
with respect to the infeasibility
for ADMM and its accelerated version,
using~\eqref{eq:ppm,rate} and~\eqref{eq:appm,rate}.

\begin{corollary}
Assume that $||\hat{\nuu}_0 + \rho\A(\x_0 - \cc) - \nuu_*|| \le R$
for some $\nuu_*\in X_*(\M_{\rho,-\A\partial f^*(-\A^*\cdot) - \cc,-\B\partial g^*(-\B^*\cdot)})$.
Alternating direction method of multipliers satisfies
\begin{align}
||\A\x_{i+1} + \B\z_i - \cc||^2 \le \left(1 - \frac{1}{i}\right)^{i-1}\frac{R^2}{\rho^2i}
\label{eq:admm,rate}
,\end{align}
and the proposed accelerated alternating direction method of multipliers satisfies
\begin{align}
||\A\x_{i+1} + \B\z_i - \cc||^2 \le \frac{R^2}{\rho^2i^2}
\label{eq:aadmm,rate}
.\end{align}
\end{corollary}

The bound~\eqref{eq:admm,rate} is 
$e$-times asymptotically
smaller than the known rate for ADMM
in~\cite[Theorem 15]{davis:16:cra},
which originated from the bound~\eqref{eq:ppm,mono}.
Finding exact bounds for the ADMM and its proposed variant
is yet left as future work.

\begin{remark}
Many existing rates for the (preconditioned) ADMM
consider the \emph{ergodic} sequences $\{\bar{\x}_i\}$ and $\{\bar{\z}_i\}$,
where $\bar{\x}_i := \frac{1}{i}\sum_{l=1}^i\x_l$ 
and $\bar{\z}_i := \frac{1}{i}\sum_{l=1}^i\z_l$
(see \eg, \cite{chambolle:11:afo,chambolle:16:ote,davis:16:cra,davis:17:fcr}).
In particular, in~\cite[Theorem 15]{davis:16:cra},
ADMM is found to satisfy
\begin{align}
||\A\bar{\x}_{i+1} + \B\bar{\z}_i - \cc||^2 \le \frac{16R^2}{\rho^2i^2}
\label{eq:ergo,feas}
,\end{align}
which is faster than the rate of the \emph{nonergodic} sequence $\{\x_i,\z_i\}$ 
of ADMM in~\eqref{eq:admm,rate}
and is comparable to the rate of the proposed accelerated ADMM in~\eqref{eq:aadmm,rate}.
One should note that the feasibility convergence of the \emph{ergodic} sequence,
as in~\eqref{eq:ergo,feas},
does not necessarily
imply the convergence of the fixed-point residual of the \emph{ergodic} sequence,
unlike~\eqref{eq:admm,rate} and~\eqref{eq:aadmm,rate} for the \emph{nonergodic} sequence.
In addition,
some numerical experiments in~\cite{chambolle:16:ote}
illustrate that the performance of the \emph{nonergodic} sequence
can be faster than that of the \emph{ergodic} sequence.
We leave further understanding the rates of the \emph{ergodic} and \emph{nonergodic} sequences
of (preconditioned) ADMM and their relationship as future work.
\end{remark}

\begin{remark}
\cite{chambolle:11:afo,chambolle:16:ote,goldstein:14:fad}
proposed accelerated variants of (preconditioned) ADMM
under some additional conditions,
while the proposed method does not
require such conditions.
\end{remark}

\begin{example}
We apply the accelerated ADMM to the problem
\begin{align}
\min_{\x\in\reals^{d_1},\z\in\reals^{d_2}}\; & \frac{1}{2}||\H\x - \bb||^2 + \gamma||\z||_1 
\label{eq:tv,constrained} \\
\text{subject to }\; & \D\x - \z = \zero, \nonumber
\end{align}
with a positive real number $\gamma$,
associated with the total-variation-regularized least-squares problem
\begin{align}
\min_{\x\in\reals^{d_1}}\; & \frac{1}{2}||\H\x - \bb||^2 + \gamma||\D\x||_1
\label{eq:tv}
,\end{align}
where 
$\H \in\reals^{p\times d_1}$, $\bb \in\reals^p$,
and
a matrix $\D \in\reals^{d_2\times d_1}$ is given as
\begin{align}
\D = 
\left[\begin{array}{cccccc}
1 & -1 & 0 & 0 & \cdots & 0 \\
0 & 1 & -1 & 0 & \cdots & 0 \\
\vdots & \ddots & \ddots & \ddots & \ddots & \vdots \\
\vdots & \ddots & 0 & 1 & -1 & 0 \\
0 & \cdots & \cdots & 0 & 1 & -1
\end{array}\right]
.\end{align}
By letting $f(\x) = \frac{1}{2}||\H\x - \bb||^2$, $g(\z) = \gamma||\z||_1$,
$\A = \D$, $\B = -\I$ and $\cc = \zero$,
we have the following accelerated ADMM method:
\begin{align}
\x_{i+1} &= \argmin{\x\in\reals^{d_1}} \left\{\frac{1}{2}||\H\x - \bb||^2 + \inprod{\hat{\nuu}_i}{\D\x - \z_i}
                + \frac{\rho}{2}||\D\x - \z_i||^2\right\} \\
		&= (\H^\top\H + \rho\D^\top\D)^{-1}(\D^\top(\rho\z_i - \hat{\nuu}_i) + \H^\top\bb) 
		\nonumber \\
\hat{\etaa}_i &= \begin{cases}
		\hat{\nuu}_i, & i=0,1, \\
		\hat{\nuu}_i + \frac{i-1}{i+1}(\hat{\nuu}_i - \hat{\nuu}_{i-1} + \rho\D(\x_{i+1} - \x_i)) 
                - \frac{i-1}{i+1}(\hat{\nuu}_{i-1} - \hat{\etaa}_{i-2} + \rho\D(\x_i - \x_{i-1})), 
		& i=2,3,\ldots
		\end{cases}
		\nonumber \\
\z_{i+1} &= \argmin{\z\in\reals^{d_2}} \left\{\gamma||\z||_1 + \inprod{\hat{\etaa}_i}{\D\x_{i+1} - \z} 
                + \frac{\rho}{2}||\D\x_{i+1} - \z||^2\right\} 
		= \Soft_{\frac{\gamma}{\rho}}\left(\D\x_{i+1} + \frac{1}{\rho}\hat{\etaa}_i\right)
		\nonumber \\
\hat{\nuu}_{i+1} &= \hat{\etaa}_i + \rho(\D\x_{i+1} - \z_{i+1}), \nonumber
\end{align}
where the soft-thresholding operator is
defined as $\Soft_{\tau}(\z) := \max\{|\z|-\tau,\zero\} \odot \sign(\z)$
with the element-wise absolute value, maximum and multiplication operators, 
$|\cdot|$, $\max\{\cdot,\cdot\}$ and $\odot$, respectively.

In the experiment, we choose $d_1 = 100$, $d_2 = 99$, $p=5$, and 
a true vector $\x_{\mathrm{true}}$ is constructed
such that a vector $\D\x_{\mathrm{true}}$ has few nonzero elements.
A matrix $\H$ is randomly generated
and a noisy vector $\bb$ is generated by adding randomly generated (noise) vector 
to $\H\x_{\mathrm{true}}$.
We choose the parameters $\gamma = 3$ and $\rho = 0.05$ in the experiment.

Figure~\ref{fig:tv}
illustrates the fixed-point residual 
of ADMM and its accelerated variants.
Interestingly, ADMM has a rate comparable to the $O(1/i^2)$ rate of the proposed method.
This does not contradict with the theory,
and we leave further investigating the worst-case rate of ADMM
under the Lipschitz continuity condition of $\nabla f$;
similar analysis but under different conditions can be found in~\cite{davis:16:cra,davis:17:fcr}.
Noticing the oscillating behavior of the proposed ADMM in Fig.~\ref{fig:tv},
we heuristically restarted the proposed method every $20$ iterations,
yielding a linear rate, 
without a strong monotonicity condition.\footnote{Since $\nabla f$ is Lipschitz continuous,
the operator $\M_1 = -\D\partial f^*(-\D^\top\cdot)$ in~\eqref{eq:M1M2}
for the problem~\eqref{eq:tv,constrained} 
is strongly monotone,
but this is insufficient to guarantee a strong monotonicity of
$\M_{\rho,\M_1,\M_2}$~\eqref{eq:Mrho} for the problem~\eqref{eq:tv,constrained}.}
Restarting has been previously found useful for a different accelerated ADMM in~\cite{goldstein:14:fad}.

\begin{figure}[h!]
\begin{center}
\includegraphics[clip,width=0.55\textwidth]{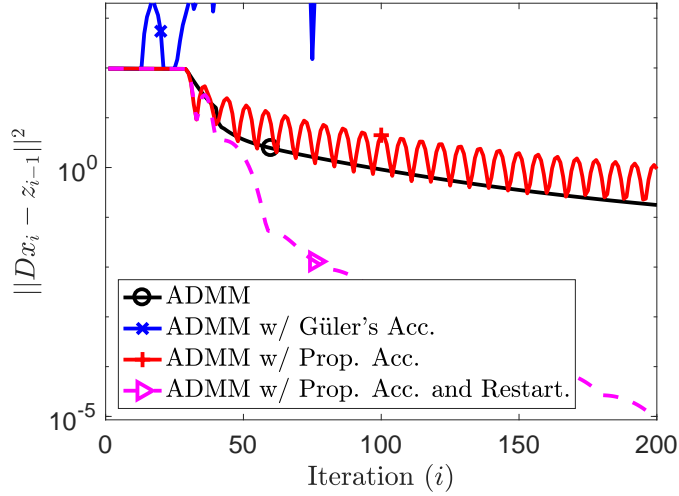}
\end{center}
\caption{Solving a total-variation-regularized least-squares problem~\eqref{eq:tv};
the fixed-point residual vs. iteration}
\label{fig:tv}
\end{figure}

\end{example}

\section{Accelerated forward method for cocoercive operators}
\label{sec:disc}

This section applies the proposed acceleration to a forward method, such as a gradient method, for cocoercive operators.
A single-valued operator $\M\;:\;\cH\to\cH$ is $\beta$-cocoercive for $\beta\in\reals_{++}$ if
\begin{align}
\inprod{\x - \y}{\M\x - \M\y} \ge \beta||\M\x - \M\y||^2\text{ for all } \x,\y\in\cH.
\label{eq:cC}
\end{align}
Let $\cC$ 
be the class of $\beta$-cocoercive operators on $\cH$.
For the $\beta$-cocoercive operator, 
the following forward method 
(that iteratively applies the forward operator $\I - \beta\M$)
is guaranteed to converge weakly to a solution
\cite[Theorem 26.14]{bauschke:11:caa}.

\fbox{
\begin{minipage}[t]{0.85\linewidth}
\vspace{-10pt}
\begin{flalign*}
&\quad \text{\bf Forward Method} & \\
&\qquad \text{Input: } \M\in\cC,\; \y_0\in\cH. & \\
&\qquad \text{For } i = 0,1,\ldots & \\
&\qquad \qquad \y_{i+1} = (\I - \beta\M)\y_i. &
\end{flalign*}
\end{minipage}
} \vspace{5pt}

An operator $\T\;:\;\cH\to\cH$ is $\lambda$-cocoercive if and only 
if it is the Yosida approximation of index $\lambda$~\cite[Proposition 23.21]{bauschke:11:caa}:
\begin{align}
\Mlam := \frac{1}{\lambda}(\I - \J_{\lambda\M})
\end{align}
of a maximally monotone operator $\M\;:\;\cH\to2^\cH$.
We thus have the following equivalence
between the resolvent (backward) operator of a maximally monotone operator $\M$
and a forward operator of the corresponding cocoercive operator $\Mlam$:
\begin{align}
\J_{\lambda\M} = (\I + \lambda\M)^{-1} = \I - \lambda\Mlam
.\end{align}
Therefore, the results on the proximal point method and its accelerated variant
for monotone operators
directly apply to the forward method and its accelerated variant
for cocoercive operators.

\section{Conclusion}
\label{sec:conc}

This paper developed an accelerated proximal point method 
for maximally monotone operators,
with respect to the fixed-point residual,
using the computer-assisted performance estimation problem approach. 
Restarting technique was further employed
under the strong monotonicity condition.
The proposed acceleration 
was applied 
to various instances of the proximal point method
such as the proximal method of multipliers,
the primal-dual hybrid gradient method,
the Douglas-Rachford splitting method,
and the alternating direction method of multipliers,
yielding accelerations both theoretically and practically.
The acceleration was also applied to
a forward method for cocoercive operators.

We leave developing accelerations
for more general or more specific classes of problems or methods
as future work,
possibly via the performance estimation problem approach;
a comprehensive understanding of 
accelerations for the alternating direction method of multipliers
with respect to various performance measures
under various conditions
are yet remain open.

\begin{acknowledgements}
The author sincerely appreciates the useful comments by the associate editor and anonymous referees.
The author also would like to thank Dr. Felix Lieder,
who brought to attention his Ph.D. thesis~\cite{lieder:18:pbm} 
and his paper~\cite{lieder:20:otc}, after the acceptance of this paper,
which optimized the step coefficients of the Krasnosel'skii-Mann iteration
for a nonexpansive operator $\T$, similarly using the PEP approach.
The form of the resulting optimized method 
differs from that of the accelerated proximal point method proposed in this paper,
but \cite{ryu:20} recently showed that they are equivalent in the sense that they generate the same sequence,
when $\T = 2\J_{\lambda\M} - \I$.
\end{acknowledgements}

\bibliographystyle{spmpsci}
\bibliography{master,mastersub}
\end{document}